%% file: ms.tex
\documentclass{article}

\usepackage{geometry}
\usepackage{bm}
\usepackage{amsfonts}
\usepackage{amsmath}
\usepackage{amssymb}
\usepackage{color}
\usepackage{verbatim}
\usepackage[colorlinks=true, allcolors=blue]{hyperref}
\usepackage{multicol}
\usepackage{enumerate}
\usepackage{color}
\usepackage{graphicx}
\usepackage{float}
\usepackage{caption}
\usepackage{cite}
\usepackage{algorithm}
\usepackage{algorithmic}
\usepackage{amsthm}

\graphicspath{{./figures/}{./exploiting_oscillations/figs/}}
\numberwithin{equation}{section}

\newtheorem{theorem}{Theorem}[section]

\newtheorem{definition}[theorem]{Definition}
\newtheorem{example}[theorem]{Example}
\newtheorem*{remark}{Remark}

\newcommand*{\supp}{\mathrm{supp}}
\newcommand*{\argmin}[1]{{\underset{#1}{\operatorname{argmin}}}}

\title{ACCELERATED VARIANCE BASED JOINT SPARSITY RECOVERY OF IMAGES FROM FOURIER DATA}
\author{Theresa Scarnati \footnote{Air Force Research Laboratory, Wright-Patterson AFB, OH, theresa.scarnati.1@us.af.mil}, Anne Gelb \footnote{Dartmouth College, Department of Mathematics, Hanover, NH, annegelb@math.dartmouth.edu}} 
\date{}

\begin{document}
\maketitle

\begin{abstract}
\input{abstract}
\end{abstract}

%

\input{introduction}

\section{Preliminaries}
\label{sec:preliminaries}
For synthetic aperture radar (SAR), magnetic resonance (MR), ultrasound and other imaging applications, data are typically modeled as the Fourier coefficients corresponding to the underlying function $f$ that we seek to recover {\cite{Gotcha,Cetin2005,Moses2004,sanders2017composite,ao2017sparse,brainweb2,LDP}}. For ease of presentation, we consider a one-dimensional piecewise smooth, 2$\pi$-periodic function $f:\mathbb{R}\rightarrow[-\pi,\pi]$. As will be demonstrated in our numerical experiments, it is straightforward to extend our technique to multiple dimensions, as well as in other intervals.  

Suppose we acquire the first $2N + 1$ noisy Fourier coefficients, $\hat{\bm f} \in \mathbb{C}^{{2N+1}}$, given elementwise as
\begin{equation}
\label{eq:four_coef}
\hat{f}_k = \frac{1}{2\pi}\int_{-\pi}^{\pi} f(x) e^{-ikx} dx + {\nu}_k, \quad k = -N,...,N,
\end{equation}
where $\bm\nu \sim {\mathcal{CN}}[\mu,\sigma]$ is complex additive Gaussian noise with mean $\mu$ and standard deviation $\sigma$.  We seek to recover ${\bm f}\in\mathbb{R}^{{N_x}}$, where each element  ${f}_j$ of $\bm{f}$ should approximate each $f(x_j), j = 1,...,N_x$, with
\begin{equation}
\label{eq:grid1D}
x_j = -\pi + \frac{2\pi(j-1)}{N_x}.
\end{equation}
{Here ${\bm x} = \{x_j\}_{j = 1}^{N_x}$ is chosen to be uniformly spaced for simplicity, and is not required for our algorithm.  {We will also choose $N_x = 2N$ since the acquired data are Fourier samples.} It is important to note that in practice (\ref{eq:four_coef}) can only be numerically modeled using the discrete Fourier transform operator, that is, 
\begin{equation}
\label{eq:fourmodel}
\mathcal{F}{\bm f} {\approx} \hat{\bm f},
\end{equation} 
where $\hat{\bm f} = \{\hat{f}_k\}_{k=-N}^N$, {$\mathcal{F}\in\mathbb{C}^{2N+1 \times N_x}$ and each element of $\mathcal{F}$ is defined as}
\begin{equation}
\label{eq:discrete_coef}
\mathcal{F}(k,j) =  {\frac{1}{N_x}e^{-ikx_j}}, \quad k = -N,...,N, \quad j = 1,...,N_x.
\end{equation} 
{As a consequence, and even without noise, a model mismatch always occurs between $f(x_j)$ and ${\bm f}_j$, $j = 1,\cdots,N_x$.} 

We now define the corresponding jump (edge) function $[f]$ as\footnote{{Depending on the context, $[f]$ can be referred to as an edge or jump function.  We use these terms interchangeably.}}  
\begin{equation}
\label{eq:jump_funct}
[f](x) = f(x^+)-f(x^-).
\end{equation}
Observe that $[f](x) = 0$ everywhere except at a jump location.  We discretize $[f](x)$ by assuming that there is at most one edge within a cell $I_j= [x_j,x_{j+1})$ for all $j = 1,...,N_x-1$. {Let us also define ${\bm g}$ as the edge vector, such that the components of ${\bm g}$ are the jump values occurring within each corresponding cell. Specifically, the components of ${\bm g}$ are given by
\begin{equation}
\label{eq:jump_funct_disc}
{\bm g}_j = \sum_{l = 1}^{N_x} [f](x_l)\delta_l(j),\quad\quad j = 1,\cdots,N_x,
\end{equation}
where $\delta_l(j) = 1$ if $l = j$ and $0$ otherwise.  Observe that since $[f](x_j) = 0$ in almost all $I_j$, (\ref{eq:jump_funct_disc}) has {\em sparse representation}.  Formally, we have} 
\begin{definition}\cite{chen2006theoretical, daubechies2010iteratively}
\label{def:sparsity}
A vector $\bm{g} \in \mathbb{R}^{N_x}$ is $s$-sparse for some $1 \leq s \leq N_x$ if
$$
||{\bm g}||_{0} =  |\supp(\bm{g})|  \leq s.
$$
\end{definition}

The sparsity of ${\bm g}$ allows us to recover an approximation $\bm{f}^*$ of $\bm{f}$ via $\ell_1$ regularized inversion techniques, also commonly referred to as compressive sensing, \cite{donoho2006compressed,candes2006robust,candes2006stable}. For example, \cite{wasserman2015image, archibald2016image, stefan2010improved} all solve the following $\ell_1$ regularization problem:
\begin{equation}
\label{eq:l1_reg}
\bm{f}^* = \argmin{\bm{q}\in\mathbb{R}^{N_x}}\left\{ \frac{1}{2}||\mathcal{F}\bm{q} - \hat{\bm f}||_2^2 + \lambda ||\mathcal{L}^m \bm{q}||_1 \right\}.
\end{equation} 
Here $\mathcal{F}$ is the discrete Fourier transform operator defined in (\ref{eq:discrete_coef}), $\lambda > 0$ is the regularization parameter, and $\mathcal{L}^m\in\mathbb{R}^{N_x\times N_x}$ is the $m$th order Polynomial Annihilation (PA) operator, \cite{archibald2005polynomial}, with
\begin{equation}
\label{eq:PA}
\mathcal{L}^m \bm{q} \approx {\bm g}. 
\end{equation}
{Of course, other sparsifying transforms such as wavelets  may also be used, but the PA operator is sufficient for our investigation.} {For convenience, a brief explanation of how the PA transform matrix is constructed is provided in  Appendix \ref{appendix:PA}.} 

\begin{remark}
\label{rem:PA_conv}
The Polynomial Annihilation (PA) operator as it appears in \cite{archibald2016image} is very closely related to High Order Total Variation (HOTV), \cite{stefan2010improved}.  There are some small differences, however, including how they are normalized and their adaptability to non-uniform grids. For consistency we exclusively use the PA operator as designed in \cite{archibald2016image} and refer readers there for additional details. For our purposes here, it is important to note that the PA {operator is  local}, and that $\mathcal{L}^m \bm{f} = \bm g + {\mathcal O}(\Delta x)^{m}$ for all components of ${\bm g}$ {\em outside} any $m+1$ length interval surrounding a true edge location.  Hence as $m$ increases, $\mathcal{L}^m \bm{f}$  is indeed a sparse vector.  Within the $m+1$ length interval surrounding a true edge, however, there are  oscillatory (non-zero) responses, meaning that as $m$ increases there are more non-zero values near the edges. This trade off has been addressed in \cite{wasserman2015image, archibald2016image}.  As will be demonstrated in what follows, the variance based joint sparsity method, by its construction, helps to mitigate the negative impact of the oscillatory responses near the edges while reinforcing the higher order convergence in smooth regions.  
\end{remark}

\section{Joint Sparsity Recovery from Multi-Measurement Vectors}
\label{sec:MMV}

\subsection{Weighted $\ell_p$ regularization}
\label{sec:weightelp}

To motivate our discussion, let us further consider (\ref{eq:l1_reg}). A typical complaint among practitioners is that when using $\ell_1$ regularization, it is difficult to choose an optimal parameter $\lambda$ without tedious hand-tuning. {First, $\lambda$ in \eqref{eq:l1_reg} is obviously application dependent since its scaling depends on the magnitude of the given measurements. The match between the fidelity and regularization terms is also application dependent.}  Moreover, even when the data source is the same,  the reconstruction may vary wildly based on the chosen value of $\lambda$.  Finally it is often difficult to choose $\lambda$ due to noise in the acquired data, incomplete information, and not having a good approximation of the sparse vector ${\bm g}$.  

Another fundamental issue is that while ${\bm g}$ is sparse, it is not zero, and the non-zero values should be preserved for a more accurate reconstruction of the underlying piecewise smooth signal. Indeed, near discontinuities the projection ${\mathcal L}^m g$ of ${\bm g}$ has $m+1$ non-zero elements in non-smooth regions of ${\bm f}$. These values should also be preserved for a more accurate reconstruction of the underlying signal.\footnote{The projection ${\mathcal L}^1g$ is often referred to as the TV or gradient domain for $m = 1$.}  This issue has been addressed in \cite{candes2008enhancing, chartrand2008iteratively, daubechies2010iteratively, liu2012adaptive} by the introduction of weighted $\ell_p$ schemes, which are designed to reduce the penalty at locations where ${\mathcal L}^m{\bm g}$ is non-zero. The weights may also vary in scale to accommodate the different magnitudes of the non-zero elements. In general, we define the weighted $\ell_p$ scheme as
\begin{equation}
\label{eq:VBJS}
\bm{f}^* = \argmin{\bm q\in\mathbb{R}^{N_x}} \left\{\frac{1}{p}||W\mathcal{L}^m \bm q ||_{p}^p + \frac{1}{2} ||\mathcal{F}\bm q - \hat{\bm{{f}}}||_2^2 \right\},
\end{equation}
where 
\begin{equation}
\label{eq:weightmatrix}
W = \text{diag}\left(w_1,w_2,...,w_{N_x}\right) = \text{diag}(\bm{w})
\end{equation}
is the diagonal matrix containing the weights at each grid point. {The weighting vector is designed to more heavily penalize the values close to zero in the sparse domain vector so that the true non-zero components are allowed to pass through.  

In its original construction, \cite{candes2008enhancing}, for any given model ${\mathcal F}{\bm f} = \hat{\bm f}$, (\ref{eq:VBJS}) is solved iteratively with $m = 1$ (TV regularization). The weighting vector ${\bm w}$ is also iterated on, with the regularization parameter $\lambda$ now accounted for in the spatially adaptive weight vector ${\bm w}$.   It is also possible to apply (\ref{eq:VBJS}) for either $p = 1$ or $p = 2$, \cite{gelb2018reducing, churchill2018edge}.  In particular $p = 2$ is as effective since well chosen weights should essentially {\em eliminate} the penalty on non-zero elements of a sparse ${\bm g}$. Using $p = 2$ is also more computationally efficient in this case.  In our investigation,  $\hat{\bm f}$ are the acquired noisy Fourier data, $\mathcal{F}$ is given by (\ref{eq:discrete_coef}), and $\mathcal{L}^m$ is the $m$th order Polynomial Annihilation operator \cite{archibald2005polynomial}.

{Unfortunately, noise in the acquired data and having only incomplete information make it difficult to obtain a good approximation to sparse vector $\mathcal{L}^m{\bm g}$}.  Consequently, the weights in (\ref{eq:VBJS}) may not be accurately chosen, and in some cases, the approximation is {\em worse} than the original approximation in (\ref{eq:l1_reg}). This is because the penalty may be either erroneously reduced or enhanced at different locations. In iterative reweighting schemes}, \cite{candes2008enhancing, chartrand2008iteratively, xie2014reweighted}, the situation is further exacerbated since the incorrect solution is subsequently reinforced. This is particularly a problem when $p = 2$, as the poor reconstruction is more ``spread out'' over the domain. 

\subsection{Using sparse multi-measurement vectors to obtain weights}
\label{sec:MMV_weights}

Clearly the choice of weighting vector ${\bm w}$ in \eqref{eq:weightmatrix} is critical to the approximation $\bm{f}^*$ of ${\bm f}$. To this end, as discussed in the introduction, for certain applications it is possible to acquire multiple measurements of the data for a particular region of interest.  This is the framework for the well known {\em multiple measurement vector} (MMV) problem \cite{ao2017sparse, chen2006theoretical, cotter2005sparse, eldar2009robust}. Ideally, having more than one data measurement should improve the overall quality of the data processing task (e.g.~image reconstruction, classification, etc.), and as we will demonstrate in what follows, it is particularly useful in choosing weighting vector ${\bm w}$.  This is especially true when any individual data vector does not offer complete information or is polluted by noise, since, as noted above, this complicates the process in {(iteratively)} selecting ${\bm w}$. As we describe below, having multiple measurements can significantly inform our choice of weights in (\ref{eq:weightmatrix}), thereby improving the approximation $\bm{f}^*$ of ${\bm f}$ in (\ref{eq:VBJS}).

\begin{remark}
\label{rem:four_data}
While ostensibly one could consider different sources of data as the different measurements for the MMV setup, for simplicity we assume that we acquire $J$ Fourier vectors, as in (\ref{eq:four_coef}), as our measurement vectors.
\end{remark} 

\begin{remark}
\label{rem:smv}
Indeed, as discussed in Section \ref{sec:introduction}, our algorithm also applies when a single measurement is obtained. In this case, the single measurement vector is {uniquely} processed $J$ times to recover $J$ {different} approximations of edge vectors ${\bm g}$. From an algorithmic development standpoint, this is equivalent to acquiring  $J$  measurements $\hat{\bm f}$. This approach was used in \cite{scarnati2018variance} to reduce speckle in synthetic aperture radar imagery and in \cite{adcock2019joint} for recovering edges from non-uniform Fourier data.  {In this regard, our method provides a super-resolution solution for the single measurement case.}
\end{remark}

Suppose for now that we can successfully acquire $J$ unique approximations of the edge vector ${\bm g}$ in (\ref{eq:jump_funct_disc}) of a $2\pi$-periodic piecewise smooth function $f$ from given Fourier data (\ref{eq:four_coef}),\footnote{This could be from $\tilde{J} \le J$ measurement vectors.} given by
\begin{equation}
\label{eq:approx_jump_funct}
\tilde{\bm g}^j \approx {\bm g}, \quad j = 1,...,J,
\end{equation}
where $\tilde{\bm g}^j\in\mathbb{R}^{N_x}$ for each $j$.  In Section \ref{sec:jump_approx} we describe two options for calculating the approximations (\ref{eq:approx_jump_funct}). Regardless of the approximation technique, because each $\tilde{\bm g}^j$ is approximating the same edge vector ${\bm g}$ corresponding to the $2\pi$-periodic and piecewise smooth $f$, we have 
\begin{equation}
\label{eq:support}
\supp(\tilde{\bm g}^1) \approx \supp(\tilde{\bm g}^2)\approx\cdots \approx \supp(\tilde{\bm g}^J).
\end{equation}
That is, each edge vector approximation (\ref{eq:approx_jump_funct}) has non-zero values in \textit{approximately} the same locations.  For a variety of reasons, including the oscillatory effects of approximating a jump function near edges, as well as the acquired data possibly having different sources of noise, or missing bandwidths of data, it is not possible to guarantee that the supports are exactly equivalent. However, oscillations and noise variance resulting from the numerical approximation of the underlying jump function can be exploited to improve signal/image recovery. Indeed in \cite{adcock2019joint,gelb2018reducing}, the variance based joint sparsity (VBJS) algorithm was designed to exploit the similar oscillatory patterns in each of the $J$ approximations. In doing so, undesirable artifacts in jump function approximations (\ref{eq:approx_jump_funct}) are reduced, leading to a more accurate recovery of $\bm{f}$ via regularized inversion techniques. {The VBJS technique is briefly described below, with more detail provided in \cite{adcock2019joint,gelb2018reducing}.
 
\subsection{Variance based joint sparsity (VBJS) weighted $\ell_p$ regularization}
\label{sec:VBJS}

As discussed previously, the weighted $\ell_p$ regularization in (\ref{eq:VBJS}) will be effective if the weighting vector ${\bm w}$ in (\ref{eq:weightmatrix}) is chosen to heavily penalize the components of ${\bm f}$ where ${\bm g}$ is supposedly zero, while not penalizing places where ${\bm g}$ has non-zero values.\footnote{{To simplify our approach, we assume ${\mathcal L}^m {\bm g} = {\bm g}$. While the algorithm can be more finely tuned to accommodate the particular construction of ${\mathcal L}^m{\bm g}$ in (\ref{eq:VBJS}), our numerical results indicate that this extra tuning is not needed.}}  Ideally, we seek
\[
{w}_i{g}_i = 
\begin{cases} 
0 , \quad &  i \in \mathcal{K} \\
c , \quad & i \not\in \mathcal{K}, 
\end{cases}
\]
where $c >> 0$ and $\mathcal{K}=\left\{i \in [1,N_x] | {g}_i \neq 0\right\}$ is the set of all indices for which the corresponding cells contain jump discontinuities. {As we can only acquire approximations of $\bm{g}$, we more realistically seek to satisfy}
\begin{equation}
\label{eq:intuition_realistic}
{w}_i\tilde{{g}}_i^j = 
\begin{cases} 
0 , \quad &  i \in \mathcal{K} \\
c , \quad & i \not\in \mathcal{K},
\end{cases}
\end{equation}
for all $j = 1,...,J$. {To achieve such weights, we first recall the definition of the {\em minmod} operator for {${a}^j \in \mathbb{R}$:
\begin{equation}
\label{eq:minmod}
\text{minmod}\left\{{a}^1,...,{a}^J\right\} = \begin{cases} 
s \cdot \min\left(|{a}^1|, ..., |{a}^J|\right), \quad & {sgn}({a}^1) = \cdots
= {sgn}({a}^J) = s,\\
0, \quad & \text{otherwise}.
\end{cases}
\end{equation}
}
{We then define ${\bm S} \in \mathbb{R}^{N_x}$ componentwise as}
\begin{equation}
\label{eq:S}
{S}_i := \text{minmod} \left\{ \tilde{\bm g}_i^1, \tilde{\bm g}_i^2, ...,\tilde{\bm g}_i^J\right\}.
\end{equation}
We note that the \textit{minmod} operator has been used similarly for exploiting oscillatory responses near jump function approximations to ``pinpoint'' an exact edge, \cite{GT06}.  Specifically, as a consequence of the Gibbs phenomenon, approximating a jump function from Fourier data inevitably leads to oscillations near an edge.  However, by employing different choices of $\sigma$ in (\ref{eq:concsum}) for the approximation  of $[f](x)$, the vectors $\tilde{\bm g}^j$, $j = 1,\cdots,J$, will yield oscillations of different {\em signs} near the edges, so that ${\bm S}$ will be zero in those locations.  This is explained more in Section \ref{sec:jump_approx} when we discuss how to approximate each $\tilde{\bm g}^j$, $j = 1,\cdots,J$, from Fourier data. Note that previous implementations of the VBJS algorithm, \cite{gelb2018reducing,scarnati2018variance}, did not exploit {the differing oscillations among the sparsity vectors in the weight vector design.} In doing so here, the VBJS reconstruction is more robust to spurious oscillations due either to noise or other sources of error. The numerical examples in Section \ref{sec:numerics} demonstrate this phenomenology. {To the best of our knowledge, minmod thresholding has never been used for regularizing optimization problems, although it is commonly applied when using flux limiters for solving numerical conservation laws \cite{leveque1992numerical}.}

The other ingredient needed to generate weights that satisfy (\ref{eq:intuition_realistic}) requires the calculation of the sample variance across the rows of the joint sparsity matrix $\mathcal{P}\in\mathbb{R}^{N_x\times J}$, given by
\begin{equation}
\label{eq:JS_matrix}
\mathcal{P} = \begin{bmatrix}
\tilde{\bm g}^1 & \tilde{\bm g}^2 & \cdots & \tilde{\bm g}^J
\end{bmatrix}. 
\end{equation}
Each element of the variance vector $\bm{v}\in\mathbb{R}^{N_x}$ is determined as
\begin{equation}
\label{eq:variance}
v_i = \frac{1}{J} \sum_{j = 1}^J (\mathcal{P}_{i,j})^2  - \left( \frac{1}{J}\sum_{j = 1}^J \mathcal{P}_{i,j}\right)^2, \quad i = 1,...,N_x.
\end{equation}
Finally, by defining $\bm{T}\in\mathbb{R}^{N_x}$ {componentwise} as 
\begin{equation}
\label{eq:T}
T_i = \frac{|S_i v_i| }{\max_i |S_iv_i|}, \quad i = 1,...,N_x,
\end{equation}
we can prescribe the elements of the weighting vector $\bm w \in \mathbb{R}^{N_x}$ as 
\begin{equation}
\label{eq:weights}
w_i = \begin{cases} 
c, \quad & T_i \geq \tau \\
1 - T_i \quad & T_i < \tau. 
\end{cases}
\end{equation}
Here, $c= \left| \left\{ \ell | T_{\ell} \geq \tau\right\} \right|$ (with $|\cdot|$ denoting cardinaity of a set) is the number of elements in $\bm T$ that fall above the given threshold $\tau$ (i.e. the number of cells $I_i$ that the algorithm says contain edge locations). We note that (\ref{eq:weights}) also allows for the {\em separation of scales} within the weights. Moreover, as will be demonstrated in our numerical examples, this separation of scales generally works better than when simply ``masking'' the jump regions (i.e.~assigning binary weights {that may be scaled by an appropriately tuned regularization parameter}).}  Because the weights are derived from the variance across the rows of the joint sparsity matrix, this technique was coined the {\em variance based joint sparsity} (VBJS) recovery method in \cite{adcock2019joint, gelb2018reducing}. 

It is evident from (\ref{eq:weights}) that a threshold is needed to determine the existence of a jump. For all experiments, we set $\tau = \mathcal{O}(1/N)$ and demonstrate the robustness of this selection. {Since (\ref{eq:weights}) still produces a separation of scales in the weighting vector, this choice for $\tau$ is robust even in low SNR environments.  Indeed, in general this separation of scales in ${\bm w}$ is critical to help alleviate the difficulties due to noise and other inaccuracies in numerically computing $\mathcal{P}$.} 

\begin{remark}
\label{rem:new_weights}
The derivation of the weights in (\ref{eq:weights}) closely follows the work in \cite{gelb2018reducing}.  There is, however, an important distinction here that leads to {both greater accuracy as well as accelerated} convergence.  Specifically, in \cite{gelb2018reducing} we first performed $J$ initial reconstructions  $\{\tilde{\bm f}^j\}_{j = 1}^J$, using (\ref{eq:l1_reg}), which was then used to recover $J$ jump function approximation via the PA transform, $\tilde{\bm g}^j = \mathcal{L}^m\tilde{\bm f}^j$, $j = 1,\cdots,J$.  We then constructed the joint sparsity matrix (\ref{eq:JS_matrix}).  {By contrast,} here we calculate $\tilde{\bm g}^j$, $j = 1,\cdots, J$, {\em directly} from the Fourier data. Hence the calculation of the variance in (\ref{eq:variance}) is now more reflective of the true variance among measurements {(Fourier data)}. Moreover, it is not necessary to solve the inverse problem to calculate the sparsity vectors in (\ref{eq:approx_jump_funct}).  Thus, our new approach is both more accurate and efficient. 
\end{remark}


\section{Recovering Jump Function Approximations}
\label{sec:jump_approx}

{We now describe the concentration factor (CF) edge detection method approach for recovering edge vectors $\{\tilde{\bm g}^j\}_{j = 1}^J$ in (\ref{eq:approx_jump_funct}) from Fourier data (\ref{eq:four_coef}) so that weights (\ref{eq:weights}) can be accurately constructed and used in (\ref{eq:VBJS}).   The CF method can also be adapted to situations where measurements may be {missing or otherwise deemed unreliable.  In some instances, this may affect entire bandwidths of data.} The modification involves iteratively solving a convex optimization problem for each measurement vector. While not as efficient, the modified CF method  provides more accurate and stable reconstructions in such circumstances.} {More information can be found in \cite{viswanathan2012iterative}.}

\subsection{The concentration factor (CF) edge detection method}
\label{sec:CFmethod}
Given the Fourier coefficients of a $2\pi$ periodic piecewise-smooth function as in (\ref{eq:four_coef}), the concentration factor (CF) edge detection method, developed in \cite{GT1},  approximates the jump function (\ref{eq:jump_funct}) as 
\begin{equation}
\label{eq:concsum}
S_N^{\sigma}[f](x) = i \sum_{|k|\leq N} \hat{f}_k {sgn}(k)\sigma\left(\frac{|k|}{N}\right)e^{ikx} \approx [f](x).
\end{equation}
Here $\sigma(\eta)$, $\eta \in (0,1]$, is a ``concentration factor'' that enables (\ref{eq:concsum}) to ``concentrate'' at the singular support of $f$. A concentration factor $\sigma$ should satisfy the following admissibility conditions (see e.g.~\cite{GT1,GT06,viswanathan2012iterative}):
\begin{enumerate}
\item $K^{\sigma}_N(x) = \sum_{k=1}^{N} \sigma(\frac{k}{N})\sin(kx)$ is odd, 
\item $\frac{\sigma(\eta)}{\eta} \in C^2(0,1)$ (the first and second derivatives are continuous), and
\item $\int_{\epsilon}^1 \frac{\sigma(\eta)}{\eta}\rightarrow - \pi$, where $\epsilon  > 0 $ is small. 
\end{enumerate} 
In a nutshell, the first condition describes how (\ref{eq:concsum}) can be written as the convolution $K^{\sigma}_N \ast f$,\footnote{Recall that the standard Fourier partial sum approximation of $f$ can be written as $D_N \ast f$  where $D_N(x) = 1+ 2\sum_{k = 1}^{\frac{N}{2}} \cos{kx}$ is the usual Dirichlet kernel.}  the second condition provides the necessary smoothness needed for $\sigma$ so that (\ref{eq:concsum}) converges, and the third condition provides proper normalization. The convergence properties of (\ref{eq:concsum}) depend on the particular choice of the admissible concentration factor $\sigma(\eta)$. {Also observe that it is always possible to choose $\sigma$ so that by construction (\ref{eq:concsum}) will be oscillatory.}\footnote{It is of course possible to choose $\sigma$ so that (\ref{eq:concsum}) is not oscillatory, e.g.~as a Gaussian function.  We will not choose $\sigma$ in this way, {however, since then using (\ref{eq:minmod}) is essentially ineffective.}}  This is important because as described previously, (\ref{eq:S}) is most effectively employed when comparing oscillatory approximations, specifically since the oscillations will have different signs. From (\ref{eq:concsum}) we see that the components of edge vector are approximated as $\tilde{\bm g}_i = S_N^{\sigma}[f](x_i)$, $i = 1,\cdots,N_x$. 

We can now consider two different cases: (i) we acquire $J$ different data sets $\hat{\bm f}^j$, $j = 1,\cdots J$, given by (\ref{eq:four_coef}) from which we use (\ref{eq:concsum}) on any admissible $\sigma$ to obtain $\tilde{\bm g}^j$, $j = 1,\cdots, J$, or  (ii) we acquire one data set $\hat{\bm f}$ from which we use $J$ different admissible concentration factors $\sigma^j$, $j = 1,\cdots, J$, to obtain $\tilde{\bm g}^j$, $j = 1,\cdots,J$.  Regardless of whichever case is under consideration, the process through which to obtain the weights in (\ref{eq:weights}) needed for (\ref{eq:VBJS}) is the same. {We will refer to the general technique ({for either case (i) or (ii)}) as the {\em concentration factor VBJS (CF VBJS)} method.  Algorithms \ref{alg:CF_VBJS_MMV} and \ref{alg:CF_VBJS_SMV} respectively describe the CF VBJS procedure for each case.  We note that it is also trivial to combine the approaches by using a different concentration factors for each measurement vector.}

\begin{algorithm}[h!]
\caption{CF VBJS from Multiple Measurement Vectors.}
\label{alg:CF_VBJS_MMV}
\begin{algorithmic}[1]
\STATE Acquire multiple measurement vectors $\hat{\bm f}^j$, $j = 1,\ldots,J,$ according to (\ref{eq:four_coef}).
\STATE Using a single concentration factor $\sigma$, calculate $J$ jump function approximations as
\begin{equation}
\label{eq:jump_approx_mmv}
\tilde{\bm g}^j = i \sum_{|k|\leq N} \hat{ f}^j_k {sgn}(k)\sigma\left(\frac{|k|}{N}\right)e^{ikx}, \quad j = 1,...,J.
\end{equation}
\STATE Use (\ref{eq:JS_matrix}) - (\ref{eq:weights}) to calculate the spatially-adaptive weighting vector. 
\STATE Select the optimal measurement vector $\hat{\bm f}=\hat{\bm f}^{j^*}$ that solves (\ref{eq:jsdata_index}). 
\STATE With a choice of $p =1$ or $p = 2$ and PA transform order $m$, solve the CF VBJS reconstruction problem (\ref{eq:VBJS}).
\end{algorithmic}
\end{algorithm}

\begin{algorithm}[h!]
\caption{CF VBJS from a Single Measurement Vector.}
\label{alg:CF_VBJS_SMV}
\begin{algorithmic}[1]
\STATE Acquire a single measurement vectors $\hat{\bm f}$ according to (\ref{eq:four_coef}).
\STATE Using multiple concentration factors $\sigma^j$ for $j = 1,...,J$, calculate $J$ jump function approximations as
\begin{equation}
\label{eq:jump_approx_smv}
\tilde{\bm g}^j = i \sum_{|k|\leq N} \hat{ f}_k {sgn}(k)\sigma^j\left(\frac{|k|}{N}\right)e^{ikx}, \quad j = 1,...,J.
\end{equation}
\STATE Use (\ref{eq:JS_matrix}) - (\ref{eq:weights}) to calculate the spatially-adaptive weighting vector.  
\STATE With a choice of $p =1$ or $p = 2$ and PA transform order $m$, solve the CF VBJS reconstruction problem as in (\ref{eq:VBJS}).
\end{algorithmic}
\end{algorithm}

\subsection{Determining the best measurement vector}\label{sec:bestmeasurement}
In the MMV situation described in case (i), we must additionally decide which measurement vector to use in our final CF VBJS reconstruction (\ref{eq:VBJS}). {While this can be accomplished in a variety of ways, for our discussion here, we simply} select $\hat{\bm f}$ to be the measurement vector that is ``most similar'' to all of the other measurement vectors. This is the same approach used in \cite{gelb2018reducing}.  Specifically, we define a distance matrix ${\mathcal D}$ with entries
\begin{equation}
\label{eq:l2_msmts} 
\mathcal{D}_{i,j} = ||{ \tilde{\bm{g}}^i - \tilde{\bm{g}}^j }||_2,  
\end{equation}
and choose the measurement vector $\hat{\bm f} = \hat{\bm f}^{j^*}$, and, if appropriate, forward operator $\mathcal{F}=\mathcal{F}^{j^*}$, corresponding to the $j^*$th index that solves
\begin{equation}
\label{eq:jsdata_index}
j^* = {\argmin{j\in [1,J]}}\hspace{2mm}\sum_{i = 1}^J\mathcal{D}_{i,j}.
\end{equation}
{We point out that (\ref{eq:jsdata_index}) differs from \cite{gelb2018reducing}(3.10) in two important ways.  First, in its original form, the distance measure (\ref{eq:l2_msmts}) was based on the individual reconstructions of the underlying image ${\bm f}^j$, $j = 1,\cdots,J$.  Since we no longer calculate these approximations as part of the VBJS reconstruction in either Algorithm \ref{alg:CF_VBJS_MMV} or \ref{alg:CF_VBJS_SMV}, it is more convenient to base the distance measure on the approximated edge vector.  Secondly, by using the summation in (\ref{eq:jsdata_index}), we take into account how each vector relates to {\em all} of the other vectors.  In \cite{gelb2018reducing}, by contrast, only two measurement vectors were compared at a time.} 

\subsection{Iterative concentration factor design}
\label{sec:ICF}

By design, and to satisfy the admissibility conditions, standard concentration factors can be described essentially as band pass filters.  This is evident in Figure \ref{fig:exp_cf}(left).  However in many remote sensing applications where Fourier data are collected, certain bandwidths of data are often lost or corrupted due to intentional jamming of systems, incidental outside interference or other sources of instrument error \cite{nguyen2012recovery,ferreira1992incomplete}. {Hence there is potential mismatch between the Fourier data that are amplified by the concentration factor and the data that are actually measured. Consequently the edge vector approximation is  poor since necessary data are missing for the reconstruction.  This is further discussed in \cite{viswanathan2012iterative}.} 

To describe the potential difficulties when certain bandwidths of Fourier data are not available, consider the following example of a {ramp (or saw tooth) function with a single discontinuity at $x = 0$:}

\begin{example}
\label{ex:missing_band}
Consider reconstructing the following {ramp} function
\begin{equation}
\label{eq:sawtooth}
r(x) = r_0(x) := \begin{cases} 
\frac{-x-\pi}{2\pi}, \quad & -\pi \leq x \leq 0 \\
\frac{\pi -x}{2\pi}, \quad & 0 < x \leq \pi.
\end{cases}
\end{equation}
The exact Fourier coefficients of (\ref{eq:sawtooth}) are
\begin{equation}
\label{eq:ramp_FT} 
\hat{r}_k = \begin{cases} 
\frac{1}{2\pi i k}, \quad & k \neq 0 \\
0, \quad & k = 0. 
\end{cases}
\end{equation}
\end{example}
{For this example} we consider the case where we have acquired measurement vectors $\hat{\bm f}_k^j = \hat{r}_{k}$, $j = 1,\cdots, 4$, but rather than having all of the coefficients $-N < k < N$, each measurement is missing Fourier coefficients within the bandwidth of $\mathbb{K}_j=\{k|10j \leq |k|\leq 10j+20\}$. To emphasize the issue of the missing bandwidth of data, we assume that there are no other sources of noise on the remaining coefficients, that is, $\nu_k = 0$ for all $k = -N,...,N$ in (\ref{eq:four_coef}). We discretize the problem with $N = 64$ and $N_x = 128$ in (\ref{eq:discrete_coef}) and set $\mathcal{F}^j(\mathbb{K}_j,i) = 0$.

We initially solve the problem described in Example \ref{ex:missing_band} using exponential concentration factors (\cite{GT06,viswanathan2012iterative}), which are defined as
\begin{equation}
\label{eq:exp_cf}
\sigma_E^j(\eta) = C\eta e^{\left(\frac{1}{\alpha_j \eta(\eta-1)}\right)}, 
\end{equation}
where $\alpha_j$ defines the order of the concentration factor for the $j$th measurement, and $C$ is a normalizing constant defined as 
\begin{equation}
\label{eq:exp_cf_norm}
C = \frac{\pi}{\int_{1/N}^{1-1/N} \exp\left(\frac{1}{\alpha_j\tau(\tau-1)}\right)d\tau}.
\end{equation}
Figure \ref{fig:exp_cf}(left) displays exponential concentration factors with order $\alpha_j = 2j$ for $j = 1,...,4$. These concentration factors are then used to approximate the jump function (\ref{eq:concsum}), yielding $\tilde{\bm g}^j$ for $j = 1,...,4$ in (\ref{eq:approx_jump_funct}). The approximations in each of the four cases are shown in Figure \ref{fig:exp_cf}(middle), where it is evident that the missing bandwidths of data greatly affect the approximation to the jump function $[f](x)$ in \eqref{eq:jump_funct}. In particular, the variance is relatively large for a larger region surrounding the jump discontinuity. {Therefore, the resulting weights in (\ref{eq:weights}) do not lead to marked improvement when using VBJS (\ref{eq:VBJS}). This can be observed} in Figure \ref{fig:exp_cf}(right) which shows CF VBJS reconstructions for $p = 1$ and $p = 2$ in (\ref{eq:VBJS}) with $m =2$. (Note that here we have set $\mathcal{F} = \mathcal{F}^{j^*}$ where $j^*$ solves (\ref{eq:jsdata_index}).) By contrast, when {\em all} of the first $2N+1$ Fourier coefficients are known, the width surrounding the edge location is {only} narrowly affected by spurious oscillations in (\ref{eq:concsum}). Figure \ref{fig:exp_cf_no_miss}(left) displays the jump function approximations (\ref{eq:approx_jump_funct}) computed from all of the first $2N+1$ Fourier coefficients using (\ref{eq:concsum}) with the same exponential concentration factors seen in Figure \ref{fig:exp_cf}(left). In this case, because the variance among approximations is only large in regions surrounding discontinuities, Figure \ref{fig:exp_cf_no_miss}(right) demonstrates that the CF VBJS technique is successful at recovering the function for both $p = 1$ and $p = 2$.

\begin{figure}[htb]
\centering
\includegraphics[width = \textwidth]{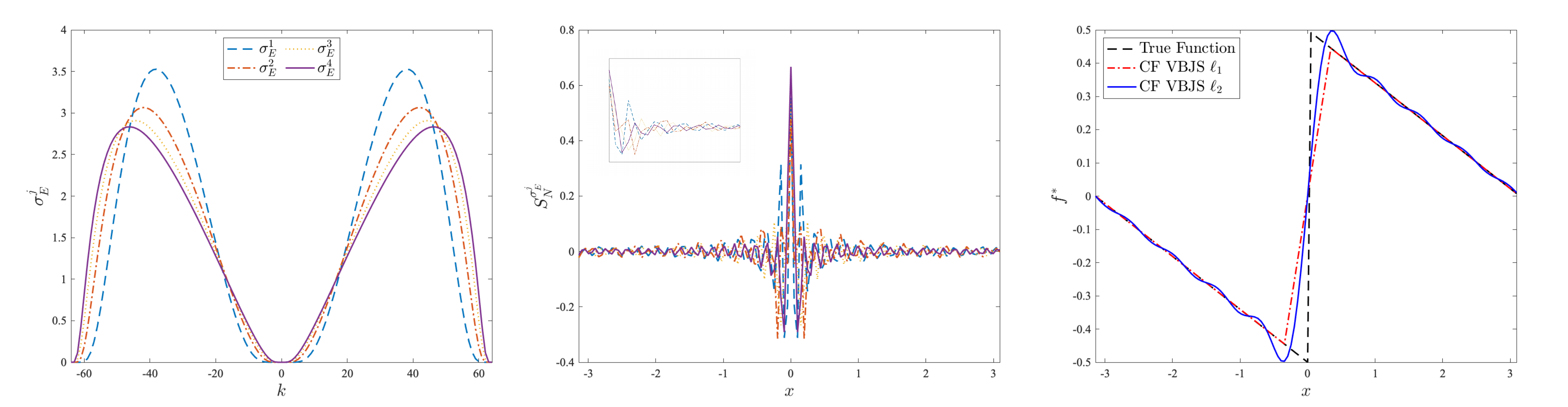}
\caption{The result of reconstructing (\ref{eq:sawtooth}) from $J=4$ measurements {of $2N+1$ Fourier samples (here $N = 64$) with missing bandwidths of data in the regions $\mathbb{K}_j=\{k|10j \leq |k|\leq 10j+20\}$},  using exponential concentration factors (\ref{eq:exp_cf}) to approximate the jump function of each measurement.  (left) The exponential concentration factors. (middle) The jump function approximations (\ref{eq:approx_jump_funct}) estimated using (\ref{eq:concsum}) {for missing bandwdiths $\mathbb{K}_j$.} {In the top left of the figure, the zoomed in view shows the approximations for $x\in[0,1]$}. (right) The CF VBJS reconstructions using $p = 1$ and $p = 2$ in (\ref{eq:VBJS}).}
\label{fig:exp_cf}
\end{figure}

\begin{figure}[htb]
\centering
\includegraphics[width = .8\textwidth]{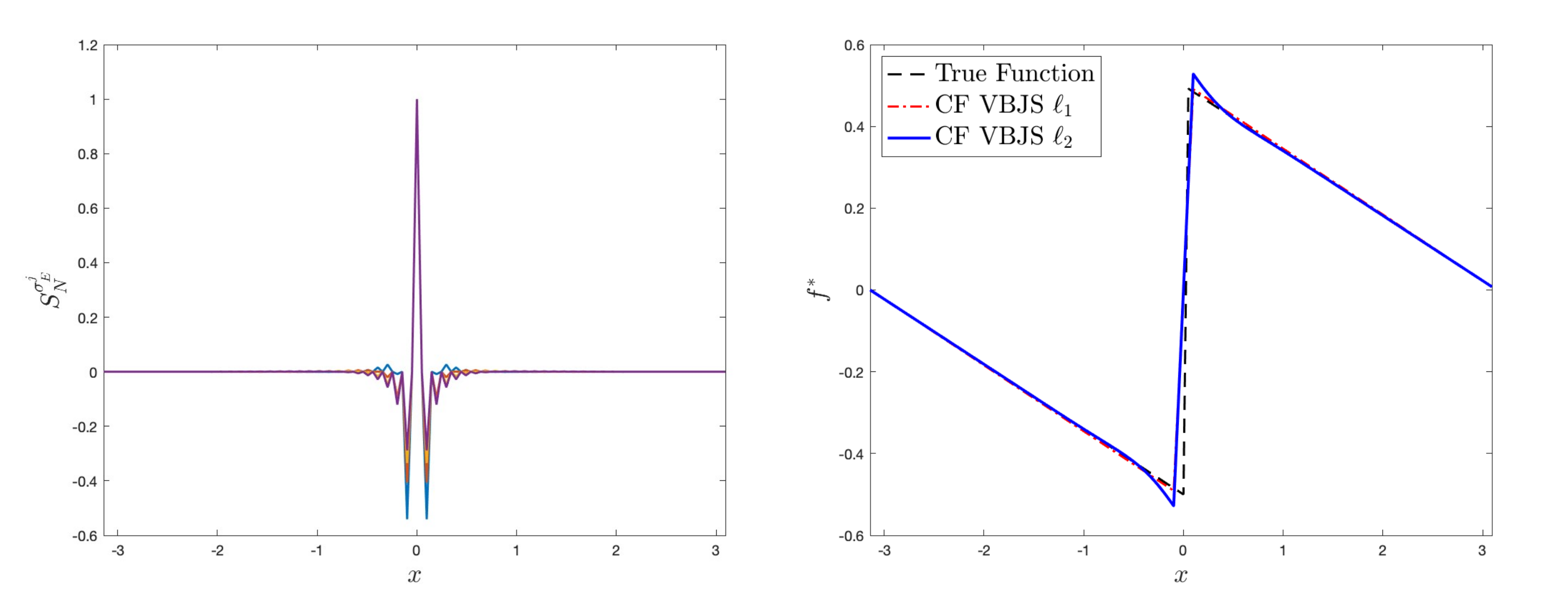}
\caption{The result of reconstructing (\ref{eq:sawtooth}) from $J=4$ measurements containing all of the first $2N+1$ Fourier coefficients (no missing bandwidths) using exponential concentration factors (\ref{eq:exp_cf}) to approximate the jump function of each measurement. {Once again $N = 64$.} (left) The jump function approximations (\ref{eq:approx_jump_funct}) estimated using (\ref{eq:concsum}) (right) The CF VBJS reconstructions using $p = 1$ and $p = 2$ in (\ref{eq:VBJS}).}
\label{fig:exp_cf_no_miss}
\end{figure}

To combat these issues, we propose {using the convex optimization} concentration factor design framework {developed in \cite{viswanathan2012iterative}} to determine a unique concentration factor {vector, ${\bm \sigma} = \{\sigma_k\}_{k=-N}^N$}, for each measurement vector {where bandwidths of Fourier data are not available}.  

{To this end, we first observe that the jump function approximation (\ref{eq:concsum}) of any periodic function $f$ defined on $[-\pi,\pi)$ for which there is a single discontinuity at $x = \xi$} can be equivalently expressed as 
\begin{equation}
\label{eq:sum_w}
S_N^{\sigma}[f](x) = \sum_{l=0}^{\infty} [f^{(l)}](\xi)W_l^{\sigma}(x-\xi), 
\end{equation}
where $\xi$ indicates {the} jump location, $[f^{(l)}]$ refers to the corresponding jump function of the $l^{th}$ derivative of $f$,  and
\begin{equation}
\label{eq:jump_Wq}
W_q^{\sigma}(x) = \frac{1}{2\pi i^q} \sum_{0<|k|\leq N} \frac{{sgn}(k)}{k^{q+1}}\sigma\left(\frac{|k|}{N}\right)e^{ikx}.
\end{equation}
{From (\ref{eq:concsum}) {and \eqref{eq:ramp_FT}} it is straightforward to show that $W_0^{\sigma}(x)$ is the jump function approximation of the ramp function (\ref{eq:sawtooth}) where the discontinuity is located at $\xi = 0$, \cite{viswanathan2012iterative}.  Similarly, $W_0^{\sigma}(x -\xi)$ is the jump function approximation of (\ref{eq:sawtooth}) for any discontinuity location $\xi \in [-\pi,\pi)$.  We note that due to the linearity of (\ref{eq:concsum}) on the Fourier data, the results that follow for a single jump discontinuity hold for multiple jumps.  Of course the global nature of the Fourier data will cause additional interfering oscillations in the neighborhoods of each jump discontinuity.}

By examining (\ref{eq:sum_w}), it is evident that to recover an accurate jump function approximation of any periodic function $f$ defined on $[-\pi,\pi)$, one must emphasize $W_0^\sigma(x-\xi)$ while suppressing the impact of $W_q^\sigma(x-\xi)$ for $q > 0$.  This then provides a mechanism for designing concentration factors that allow for a more accurate approximation (\ref{eq:sum_w}) of (\ref{eq:jump_funct}).   Specifically we seek concentration factor ${\bm \sigma}$ that yield the best approximation to (translating from $x = \xi$ to $x = 0$):
\begin{equation}
\label{eq:jump_W}
{[r_0](x)  \approx W_0^{\sigma}(x) = \frac{1}{2\pi} \sum_{0<|k|\leq N} \frac{{sgn}(k)}{k}\sigma\left(\frac{|k|}{N}\right)e^{ikx}.}
\end{equation}
Noting that the jump function  of (\ref{eq:sawtooth}) is given by
\begin{equation}
\label{eq:jump_sawtooth}
[r_0](x) = \begin{cases} 
1, \quad & x = 0 \\
0, \quad & x \neq 0,
\end{cases}
\end{equation}   
{establishes the following optimization problem for determining ${\bm \sigma}= \{\sigma_k\}_{k = -N}^N$:}
\begin{subequations}\label{eq:ICF}
\begin{align}
\label{eq:it_cf_of}
\min_\sigma \hspace{2mm} & || W_0^{\sigma}(x) ||_1 &&\\
\label{eq:jump_height}
\text{subject to} \hspace{2mm} &|W_0^\sigma(0) -1 | &&\leq \delta_1 \\
\label{eq:zero_conv}
&|W_0^\sigma(x)|_{|x|\geq \delta_2} &&\leq \delta_3 \\
\label{eq:missing_bands}
&\sigma(\mathbb{K}) &&\leq \delta_4.
\end{align}
\end{subequations}
{In our numerical examples we chose $\delta_1 = \delta_3 =10^{-3}$, $\delta_2 = .35$, and $\delta_4 = 10^{-6}$, {with $\mathbb{K}$ defining the (known) bandwidth of the acquired Fourier samples.  It is possible, of course, to choose other values for each $\delta_j$, $j = 1,\dots, 4$}.  A similar optimization problem was developed  in \cite{viswanathan2012iterative}, but was modified here to ensure convexity.} The construction of (\ref{eq:ICF}) can be explained as follows: The objective function (\ref{eq:it_cf_of}) favors a concentration factor $\sigma$ that induces sparsity in the approximation (\ref{eq:jump_W}) of (\ref{eq:jump_sawtooth}). The first constraint (\ref{eq:jump_height}) encourages $W_0^\sigma$ to have an accurate jump height at the origin. We do not require the approximation to be exact, (i.e. $W_0^\sigma(0) =1$) as in \cite{viswanathan2012iterative}, because we are only concerned with recovering accurate jump locations and not accurate jump heights. The constraint (\ref{eq:zero_conv}) enforces $W_0^\sigma$ to be close to zero away from the origin, removing unwanted oscillations away from jump locations, and the constraint (\ref{eq:missing_bands}) suppresses unwanted oscillations due to missing bandwidths of data. Note that {this process returns a concentration factor {\em vector}, ${\bm \sigma}$, which no longer satisfy the admissibility conditions.} However, as our results will demonstrate in Section \ref{sec:numerics}, this does not affect the accuracy of our reconstructions.  In sequel we refer to this concentration factor design method as the iterative concentration factor (ICF) method.

To demonstrate the utility of the ICF technique we again carry out the experiment defined in Example \ref{ex:missing_band}. The results are displayed in Figure \ref{fig:it_cf}. Figure \ref{fig:it_cf}(left) shows the concentration factors determined from (\ref{eq:ICF}) for each of the four measurements. Figure \ref{fig:it_cf}(middle) displays the resulting jump function approximations (\ref{eq:approx_jump_funct}) estimated using (\ref{eq:concsum}) for each $\mathbb{K}_j$, $j = 1,...,4$. Notice that {due to the constraints (\ref{eq:zero_conv}) and (\ref{eq:missing_bands})}, oscillations away from discontinuities are dramatically reduced. The remaining strong oscillations are also more localized, and consequently, so are the large variance values in (\ref{eq:variance}).  As a result, the weights in (\ref{eq:weights}) provide a more accurate penalization for  Algorithm \ref{alg:CF_VBJS_MMV}.  This is demonstrated in Figure \ref{fig:it_cf}(right),  where Algorithm \ref{alg:CF_VBJS_MMV} is used to reconstruct (\ref{eq:sawtooth}) via the ICF VBJS technique for both $p=1$ and $p =2$. 

\begin{remark}
\label{rem:HOICF}
{We note that the results in Figure \ref{fig:it_cf} represent the ideal case, since the ICF is designed using Fourier coefficients of the  ramp function.   More options are considered in \cite{viswanathan2012iterative} to include suppressing $W_q^\sigma$, $q > 0$, e.g. so that $||W^\sigma_q|| \le \delta_q$. This essentially provides a higher order reconstruction of $[f](x)$ {\em away} from jump discontinuities when $f$ has more variation between edges.  However, for our purposes, which is to use each approximation of $[f](x)$ in the construction of (\ref{eq:S}) and subsequently (\ref{eq:variance}) and (\ref{eq:weights}), such additional constraints are not needed.}
\end{remark}

\begin{figure}[htb]
\centering
\includegraphics[width = \textwidth]{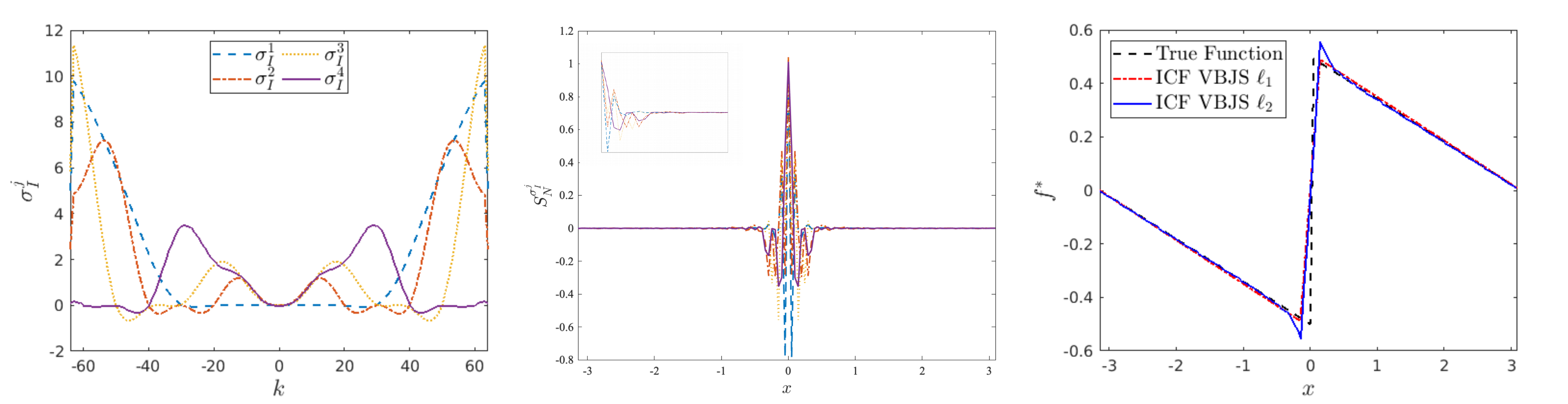}
\caption{The result of reconstructing (\ref{eq:sawtooth}) from $J=4$ measurements with missing bandwidths of data using the iterative concentration factor design technique, outlined in (\ref{eq:it_cf_of})-(\ref{eq:missing_bands}), to approximate the jump function of each measurement. (left) The iterative concentration factors. (middle) The jump function approximations (\ref{eq:approx_jump_funct}) estimated using (\ref{eq:concsum}). {In the top left of the figure, the zoomed in view shows the approximations for $x\in[0,1]$}. (right) The CF VBJS reconstructions using $p = 1$ and $p = 2$ in (\ref{eq:VBJS}).}
\label{fig:it_cf}
\end{figure}

To compare the accuracy of reconstructing (\ref{eq:sawtooth}) using exponential CF VBJS (as seen in Figure \ref{fig:exp_cf}(right)) and ICF VBJS (as seen in Figure \ref{fig:it_cf}(right)), we calculate the relative reconstruction error as
\begin{equation}
\label{eq:rel_error}
E_{rel} = \frac{||\bm{f}^* - \bm{f}||_2}{||\bm{f}||_2},
\end{equation}
and the absolute error at a grid point neighboring a discontinuity as
\begin{equation}
\label{eq:abs_error}
E_{abs} = | \bm{f}^*(x_*) - \bm{f}(x_*)|,
\end{equation}
{where $\bm{f}^*$ is calculated using (\ref{eq:VBJS}) with weights chosen via the appropriate edge approximations and $\bm{f}$ is the underlying piecewise smooth function we wish to recover on the grid defined by (\ref{eq:grid1D}).}  For each experiment, we calculate (\ref{eq:rel_error}) across {the whole domain, as well as only in regions where $f$ is at least $\mathcal{C}^1$.}  (Note that $f$ defined by (\ref{eq:sawtooth}) is in $\mathcal{C}^1$ for all $1\leq |x|\leq\pi$.   
Table \ref{tab:it_vs_exp_compare} summarizes this accuracy comparison when using $p = 1$ and $p = 2$ in the VBJS recovery (\ref{eq:VBJS}). Clearly, using the ICFs  designed in (\ref{eq:it_cf_of})-(\ref{eq:missing_bands}) allows for more accurate reconstructions when measurements are missing bandwidths of data. This is especially true in areas near discontinuities. These results will be explored further in Section \ref{sec:numerics}. 

\begin{table}[htb]
\begin{center}
\begin{tabular}{|r|c|c|c|}
\hline
 		 & Relative Error  & Relative Error & Absolute Error \\ 
 		 & (total) &  (in smooth regions)  &  ({at $x_*   \approx -0.1$})  \\ \hline
CF VBJS $\ell_1$ & 0.3325 & 0.0173  & 0.2859 \\ \hline
ICF VBJS $\ell_1$ & 0.2215 & 0.0010  & 0.0082  \\ \hline 
CF VBJS $\ell_2$ & 0.2680 & 0.0310 & 0.1617  \\  \hline
ICF VBJS $\ell_2$  & 0.2176 & 0.0036  & 0.0371 \\  \hline
\end{tabular}
\end{center}
\caption{The error resulting from estimating (\ref{eq:sawtooth}) using the iterative concentration factor method and the exponential concentration factor method in the experiment outlined in Example \ref{ex:missing_band}. Recall that here there are $J = 4$ measurements, each with different missing bands of Fourier data, and $2N = N_x = 128$. We calculate both the relative error in the entire approximation and only in smooth regions. Also, calculated is the absolute error (\ref{eq:abs_error}) near the discontinuity, {at $x_* \approx -0.1$}.}
\label{tab:it_vs_exp_compare}
\end{table}

\section{Numerical Experiments}
\label{sec:numerics}

\input{numerics}

\section{Concluding remarks} 
\label{sec:conclusions}
\input{conclusions}

\section*{Acknowledgments} 
Theresa Scarnati's work is supported in part by AFOSR LRIR \# 18RYCOR011. Anne Gelb's work is partially supported by the grants NSF-DMS 1502640, NSF-DMS 1912685, and AFOSR FA9550-18-1-0316. Approved for public release. PA Approval \#:[88AWB-2019-2699]. 

\appendix
\section{Polynomial Annihilation (PA) Matrix Construction}
\input{appendixPA}

\bibliographystyle{plain}
\bibliography{ms}

\end{document}

%% file: abstract.tex
Several problems in imaging acquire multiple measurement vectors (MMVs) of Fourier samples for the same underlying scene.  Image recovery techniques from MMVs aim to exploit the joint sparsity across the measurements in the sparse domain.  This is typically accomplished by extending the use of $\ell_1$ regularization of the sparse domain in the single measurement vector (SMV) case to using $\ell_{2,1}$ regularization so that the ``jointness'' can be accounted for.  Although effective, the approach is inherently coupled and therefore computationally inefficient.  The method also does not consider current approaches in the SMV case that use spatially varying weighted $\ell_1$ regularization term.  The recently introduced variance based joint sparsity (VBJS) recovery method uses the variance across the measurements in the sparse domain to produce a weighted MMV method that is more accurate and more efficient than  the standard $\ell_{2,1}$ approach.  The efficiency is due to the decoupling of the measurement vectors, with the increased accuracy resulting from the spatially varying weight.   This paper introduces the accelerated VBJS which reduces computational cost even further by eliminating the requirement to first approximate the underlying image in order to construct the weights.  Eliminating this preprocessing step  moreover reduces the amount of information lost from the data, so that our method is more accurate.  Numerical examples provided in the paper verify these benefits.

%% file: introduction.tex
\section{Introduction}\label{sec:introduction}

There are many applications for which multiple (indirect) measurement vectors (MMVs) of data are acquired that represent the same underlying scene, e.g., a signal or image.   The scene is typically recovered using each {\em single} measurement vector (SMV) separately, that is without exploiting any  measurement redundancies or inter-signal correlations between measurements.  In the case of noisy or incomplete data sources, compressed sensing (CS) algorithms employing $\ell_1$ regularization  are often used so that the results are optimally sparse in some domain, e.g., the gradient or wavelet domain {\cite{donoho2006compressed,candes2006robust}}. Some processing may follow after these individual reconstructions to infer information about the underlying scene, but there is always loss of information due to the individual processing of the indirect data, \cite{wasserman2015image,archibald2016image,churchill2018edge}.

More recently a number of methods have been developed to exploit the \textit{joint sparsity} across the MMVs of the same scene via $\ell_{2,1}$ regularization. One can view $\ell_{2,1}$ regularization as an extension of the usual $\ell_1$ regularization approach.  In particular, $\ell_{2,1}$ minimization enforces the overlapping sparse support of functions in a projected domain, and reconsrtucts multiple approximations of the underlying function. This technique was developed and thoroughly analyzed throughout \cite{eldar2009robust, eldar2010average, deng2012group, zhang2009user, chen2006theoretical, singh2016weighted, zheng2012subspace} and references therein, but results depend on careful hand tuning of the parameters. {The methods cited above were designed to produce a more accurate collection of estimates but typically do not} produce a single, representative reconstruction of the underlying scene.

There have also been several developments for CS  algorithms {in the SMV case} using {spatially adaptive} weighted $\ell_1$ regularization.  The idea here is to enforce more regularization in regions where the underlying signal or image is presumably zero (without value) in the sparsity domain, and by contrast less penalty at locations in the sparsity domain corresponding to non-zero values.  These algorithms thus ostensibly improve the {accuracy and} robustness of classic $\ell_1$ regularization techniques by eliminating the need to hand tune sensitive regularization parameters.   Most of these algorithms  solve a sequence of weighted $\ell_1$ minimization problems, with weights iteratively updated at each step \cite{candes2008enhancing, chartrand2008iteratively, liu2012adaptive, chan2000high, xie2014reweighted,  yang2016enhancing}.  In some cases, \cite{churchill2018edge, daubechies2010iteratively, gorodnitsky1997sparse, wipf2010iterative}, a well chosen weight can eliminate the need to use $\ell_1$ regularization, and the optimization can be performed using the much more computationally efficient $\ell_2$ regularization.   Due to noise in the data and an inaccurate approximation  of the solution in the sparse domain, the weighted $\ell_{{p}}$, $p = 1,2$, approach is not always effective, however, and in some cases may even yield {\em worse} results than the traditional (uniformly weighted) approach.  This issue is discussed in more detail in Section \ref{sec:MMV}.

A few investigations have combined the ideas of weighted regularization with those that exploit joint sparsity for MMV, \cite{singh2016weighted, zheng2012subspace, fu2016hyperspectral, gelb2018reducing, adcock2019joint}. Noteably, the variance based joint sparsity (VBJS) technique was developed in \cite{adcock2019joint,gelb2018reducing} to use the {\em variance} between measurements in the sparsity domain to determine a spatially adaptive weight to use in the aforementioned weighted $\ell_p$, $p = 1,2$ regularization.  VBJS proved to be robust with respect to noise, and in \cite{gelb2018reducing} it was also shown to be effective even when some measurements from multiple data sources contained misleading or incorrect information. In \cite{scarnati2018variance}, VBJS was adapted for synthetic aperture radar (SAR) image formation. 

In its original form, the VBJS method used initial reconstructions from each SMV (typically using the CS framework)  {\em before}  {calculating the variance between measurements in the sparse domain.  In other words, the initial CS approximations were simply improved upon by exploiting the joint sparsity properties.} This initial processing was done regardless of how the data were acquired.  Though accurate and generalizable to any linear or non-linear forward model, the method becomes increasingly computationally complex as the number of measurement vectors increases -- a problem we describe in more detail in Section \ref{sec:numerics}.  Moreover, such {initial} processing causes information loss, especially in low {resolution and low} signal-to-noise (SNR) environments.  {Thus we are motivated to improve both the performance and efficiency of VBJS by eliminating the need to pre-process the reconstruction given certain types of data acquisition methods.} 

As will be demonstrated in this investigation, in applications where data are acquired as Fourier samples, the joint sparsity properties using VBJS can be exploited {\em without} having to first perform multiple reconstructions of the underlying image or signal.  {In so doing} we are introducing an accelerated VBJS method that both improves the efficiency and accuracy of the VBJS method  given Fourier measurements.  {To this end, we note that there are many examples for which Fourier data are acquired, including synthetic aperture radar (SAR), magnetic resonance (MR) and sonar imaging.  We also note that other data collection techniques may similarly treated, and that in this regard, Fourier data serve only as a measurement prototype.}

Our specific approach is to determine an accurate approximation of the {edges (internal boundaries)} of the underlying signal or image directly from the acquired Fourier data using the concentration factor (CF) method \cite{GT1, gelb2000detection}.  {Since edges are assumed to be sparse, in this way} we obtain an accurate projection of the unknown {signal or image}  into an appropriate sparsity domain {\em without} having to solve the inverse problem for each SMV. {As a consequence} we are able exploit joint sparsity properties in the MMV framework without reconstructing the unknown signal or image, thereby increasing the efficiency and accuracy of the VBJS algorithm proposed in \cite{adcock2019joint,scarnati2018variance, gelb2018reducing}. We then proceed as before, by determining the adaptive weights based on the variance between the multiple data sets in the sparsity domain.  

This framework also provides a method for a more robust weighting strategy in the {\em single} measurement vector (SMV) case.  In particular, the CF method can process the SMV {sparse domain} projection in {\em multiple} ways.  We can then, once again, proceed as in the MMV case.  Hence we effectively design a robust non-iterative weighted $\ell_1$ regularization method for a single source of Fourier data.  {In some sense, our new approach can be regarded as a way to introduce super-resolution into the single measurement case.} 

We will also discuss two approaches for incorporating the CF edge detection method into the VBJS algorithm.  First, we use the analytically defined concentration factors originally introduced in \cite{GT1,gelb2000detection}.  The second technique adapts the iteratively designed concentration factor approach developed in \cite{viswanathan2012iterative}.  This is necessary when bands of Fourier data are missing or otherwise deemed unsuitable.  Both approaches are applicable for either the SMV or MMV case.  

The rest of this paper is organized as follows. In Section \ref{sec:preliminaries} we provide a brief description of the Fourier data model setup and define the classic $\ell_1$ regularized inversion technique from a single measurement vector. Section \ref{sec:MMV} explains how to exploit joint sparsity from a SMV or MMVs to develop spatially adaptive weighting vectors for $\ell_p$ regularized inversion. Section \ref{sec:jump_approx} introduces the CF adapted VBJS method and considers both the analytical and iterative approaches.   Numerical experiments conducted in  Section \ref{sec:numerics} test the stability and convergence properties of our new algorithm, and compares it to  previously defined methods. We evaluate the accuracy of our method in low and high signal-to-noise environments and further demonstrate its robustness when measurements are missing bands of Fourier data. Concluding remarks are provided in Section \ref{sec:conclusions}.

%% file: numerics.tex
We now explore the accuracy and efficiency of the {(I)}CF VBJS methods described in Algorithms \ref{alg:CF_VBJS_MMV}, \ref{alg:CF_VBJS_SMV} and (\ref{eq:it_cf_of})-(\ref{eq:missing_bands}). Specifically we consider recovering piecewise continuous functions {from (i) the first $2N+1$ noisy Fourier coefficients given in (\ref{eq:four_coef}) and (ii) when some bandwidths of that data are missing. In the latter case, we further analyze the iterative concentration factor (ICF) design described in (\ref{eq:it_cf_of})-(\ref{eq:missing_bands}).

{We will compare our new algorithm to the  VBJS method in Algorithm \ref{alg:VBJS_orig}. We first observe that Algorithm \ref{alg:VBJS_orig} is comparable to the original version in \cite{adcock2019joint,gelb2018reducing} in that it first reconstructs $J$ individual approximations of the underlying function (STEP 2).  The  accuracy in Algorithm \ref{alg:VBJS_orig} is enhanced since the $minmod$ construction in (\ref{eq:S}) is used for determining the weights in (\ref{eq:weights}).  Note also that the number of measurement vectors, ${\tilde J} \ge 1$, does not have to equal  $J$, as it is possible to construct $J$ unique approximations to the underlying function using (\ref{eq:l1_reg}).  Finally, critical to each algorithm is that the weight calculation still yields the critical separation of scales for the regularization term.}
\begin{algorithm}[h!]
\caption{{Enhanced} VBJS Algorithm \cite{adcock2019joint,gelb2018reducing}}
\label{alg:VBJS_orig}
\begin{algorithmic}[1]
\STATE Acquire measurement vectors $\hat{\bm f}^j$, $j = 1,\ldots,\tilde{J},$ according to (\ref{eq:four_coef}).
\STATE Use (\ref{eq:l1_reg}) to construct $J$ approximations ${\bm f}^*$.
\STATE Determine $J$ approximations of the edge vector, (\ref{eq:approx_jump_funct}) using the polynomial annihilation (PA) approximation in (\ref{eq:EdgeDetector}).  
\STATE Use (\ref{eq:JS_matrix}) - (\ref{eq:weights}) to calculate the spatially-adaptive weighting vector.
\STATE Select the optimal measurement vector $\hat{\bm f}=\hat{\bm f}^{j^*}$ that solves (\ref{eq:jsdata_index}).
\STATE With a choice of $p =1$ or $p = 2$ and PA transform order $m$, solve the {weighted optimization} problem (\ref{eq:VBJS}).
\end{algorithmic}
\end{algorithm}

\subsection{Efficiency comparison of the VBJS methods}
\label{sec:efficiencytest}


We first test the efficiency of the CF VBJS methods. Let us consider recovering the unit ramp function (\ref{eq:sawtooth}) from its first $2N+1$ Fourier coefficients, given in (\ref{eq:ramp_FT}). Let us further assume that no noise is added to the data, i.e.~$\nu_k = 0$ for all $k$ in (\ref{eq:four_coef}). We compare the speed of recovering the unknown function (\ref{eq:sawtooth}) using (i) the  VBJS method given in  Algorithm \ref{alg:VBJS_orig}, (ii) the CF method in (\ref{eq:concsum}) using standard concentration factors in Algorithm \ref{alg:CF_VBJS_SMV}, and (iii) the CF method where the concentration factors are determined using the ICF design, (\ref{eq:ICF}). In each experiment we choose $N_x = 2N$ to be the corresponding spatial resolution.  {Since the data are noise-free, we use $J$ different concentration factors in (\ref{eq:exp_cf}) and (\ref{eq:ICF}), respectively for the CF and ICF versions, to construct the $J$ edge vector measurements in (\ref{eq:jump_approx_smv}) for  Algorithm \ref{alg:CF_VBJS_SMV}.}  {For the ICF case we  choose the same $\delta_j$, $j = 1,\cdots,4$, in (\ref{eq:ICF}) as before. While it is of course possible to choose alternative values for these parameters to obtain different concentration factors from the given $2N+1$ Fourier coefficients, as this test is for the sole purpose of measuring efficiency  we simply modify the missing bandwidths of data, i.e.~${\mathbb{K}}$ in (\ref{eq:missing_bands}), to obtain each MMV.}  

\begin{figure}[htb]
\centering 
\includegraphics[width = \textwidth]{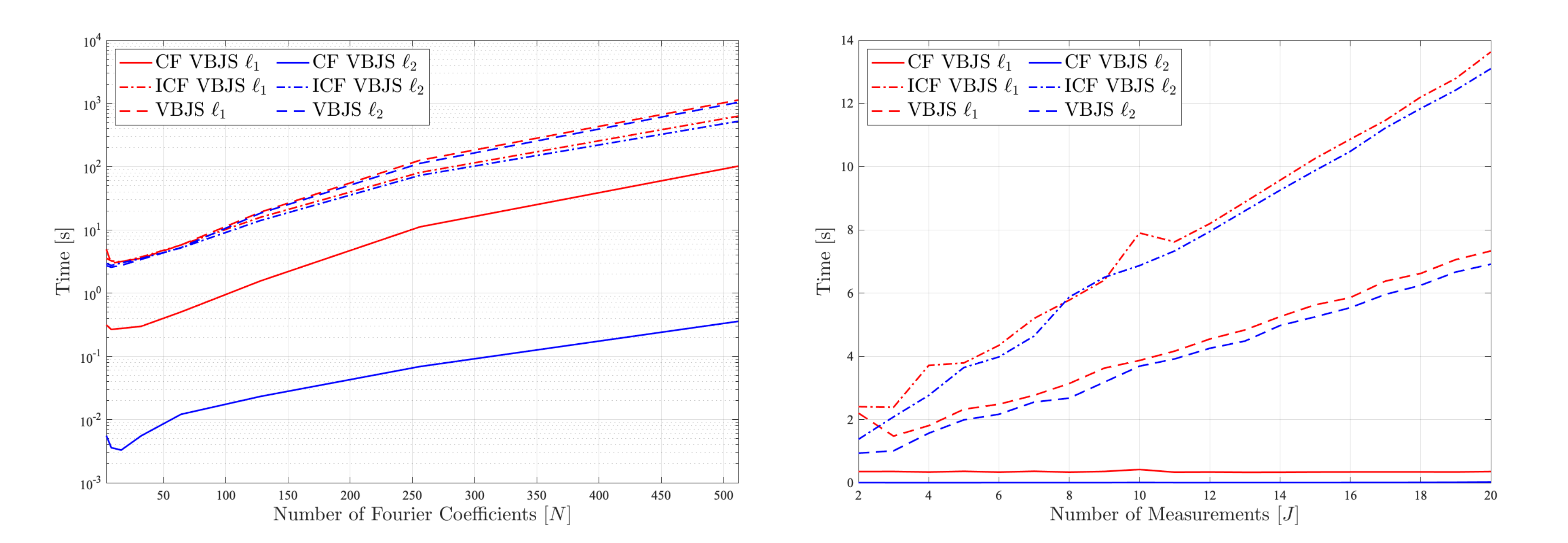}
\caption{Comparison of the VBJS recovery times (in seconds) when approximations (\ref{eq:approx_jump_funct}) of (\ref{eq:jump_sawtooth}) are found using the method posed in \cite{scarnati2018variance}, the CF method in Algorithm \ref{alg:CF_VBJS_SMV},  and the ICF design technique described in (\ref{eq:it_cf_of})-(\ref{eq:missing_bands}). In each case we acquire the first $2N+1$ noise-free Fourier coefficients. Since the data are noise-free, we use $J$ different concentration factors in (\ref{eq:exp_cf}) to construct the $J$ edge vector measurements in (\ref{eq:jump_approx_smv}) for  Algorithm \ref{alg:CF_VBJS_SMV}.  (left) $2N = N_x =2^n$, with $n = 3,4,\dots,10,$ and  $J =10$ measurements. (right) $2N = N_x = 128$ with varying number of measurements $J =2,3,\dots,20$. }
\label{fig:time_analysis}
\end{figure}

{The main distinction between the VBJS method provided in \cite{gelb2018reducing,adcock2019joint} (and enhanced by the $minmod$ step in Algorithm \ref{alg:VBJS_orig}), and our new CF VBJS technique in Algorithms \ref{alg:CF_VBJS_MMV} and \ref{alg:CF_VBJS_SMV}, is that we now {\em directly} calculate}  the $J$ edge vectors {\em without} first approximating the underlying function (in this case (\ref{eq:sawtooth})). The overall computational complexity is therefore greatly reduced.  We note that using the ICF design approach (as will be needed in the missing bandwidth case) will once again require solving $J$ optimization problems {in \eqref{eq:ICF}} to determine the concentration factors needed for (\ref{eq:jump_approx_smv}).
Figure \ref{fig:time_analysis} compares the efficiency of VBJS method in {Algorithm \ref{alg:VBJS_orig}, to our new version in Algorithm \ref{alg:CF_VBJS_SMV} with $p = 1,2$  for the regularization norm.} The time of computation (in seconds) was compared for (left) $2N = N_x =2^n$ for $n = 3,4,\dots,10$ with $J = 10$ and (right) $J = 2,3,\dots,20$ with $2N = N_x = 128$.    
The total time of each reconstruction method is measured in seconds. Each experiment was conducted using MATLAB R2018b on a single core of the same computer. 
Observe the considerable efficiency gain when using the CF VBJS technique of Algorithm \ref{alg:CF_VBJS_SMV}.  Moreover, the efficiency of the non-iterative CF VBJS technique does not strongly depend on the number of measurements acquired. As expected, using $p = 2$ is also more efficient since the corresponding objective function becomes differentiable and a closed form solution exists.


\subsection{Signal recovery from noisy Fourier data}

To further validate the CF VBJS recovery technique, we consider two experiments which aim to recover a function from noisy data. In the first experiment, we fix the noise level and test the accuracy of the method {given a varying amount of Fourier samples}. In the second experiment we fix the {the number of Fourier samples given}  and test the accuracy of the method for varying noise levels on the data. By conducting these two experiments, we are able to analyze both the convergence and stability of our technique in the presence of noise. To control the amount of noise added to the data in each case, we use the signal-to-noise ratio (SNR), defined as
\begin{equation}
\label{eq:SNR}
\text{SNR}_{\text{dB}} = 10 \log_{10} \left(\frac{E[|\hat{\bm f}|]}{\varsigma}\right).
\end{equation} 
Here $\varsigma$ is the standard deviation of the additive noise applied to the data (\ref{eq:four_coef}) and $E[|\hat{\bm f}|]$ denotes the expected value of the magnitude of the Fourier data. 

For each experiment we compare the result of the CF VBJS reconstruction to two similar techniques: (i) {the VBJS method in Algorithm \ref{alg:VBJS_orig},}  and (ii) a masking technique proposed in  \cite{churchill2019image}. In the masking technique, instead of using the weighting vector (\ref{eq:weights}) in the final VBJS reconstruction (\ref{eq:VBJS}), a binary mask is placed over grid points that fall within edge regions. Specifically, the mask is defined as
\begin{equation}
\label{eq:mask}
m_i = \begin{cases}
1, \quad & {w_i} \geq \tilde\tau \\
0, \quad & {w_i} < \tilde\tau,
\end{cases}
\end{equation} 
for {$w_i$} defined in (\ref{eq:weights}), and we set $W = \text{diag}(\bm{m})$ to solve (\ref{eq:VBJS}). {For the experiments in this paper, we set $\tilde\tau = 1$. } 

\subsubsection{Convergence experiment}
\label{ex:conv_rates}
In this first experiment with noisy data, we compare the convergence rates of the CF-VBJS method in Algorithms \ref{alg:CF_VBJS_MMV} and \ref{alg:CF_VBJS_SMV} to the VBJS technique in Algorithm \ref{alg:VBJS_orig} and the masking approach in (\ref{eq:mask}). In the multi-measurement case  we acquire $J = 10$ measurements of (\ref{eq:ramp_FT}) according to (\ref{eq:four_coef}) by independently sampling the noise vector for each measurement. That is, $\hat{f}_k^j = \hat{r}_k + \nu^j_k$ for $j=1,...,10$, with each $\bm{\nu}^j = \{\nu^j_k\}_{k=-N}^N$ sampled from a complex Gaussian distribution with mean $\mu = 0$ and standard deviation $\sigma$ such that $\text{SNR} =5$ dB. We then find the 10 edge vector approximations in (\ref{eq:jump_approx_mmv}) using the $8$th order exponential CF in (\ref{eq:exp_cf}). By contrast, in the single measurement case we are given one measurement of (\ref{eq:ramp_FT}) according to (\ref{eq:four_coef}) as $\hat{f}_k = \hat{r}_k + \nu_k$, with $\bm{\nu}= \{\nu_k\}_{k=-N}^N$ sampled from a complex Gaussian distribution with mean $\mu = 0$ and standard deviation $\sigma$ such that $\text{SNR} =5$ dB. To proceed we determine $J = 10$ edge vectors in (\ref{eq:jump_approx_smv}) using $10$ unique exponential CFs defined by (\ref{eq:exp_cf}) with order $\alpha = 2j$ for $j =1,...,10$. 

\begin{figure}[htb]
\centering
\includegraphics[width = \textwidth]{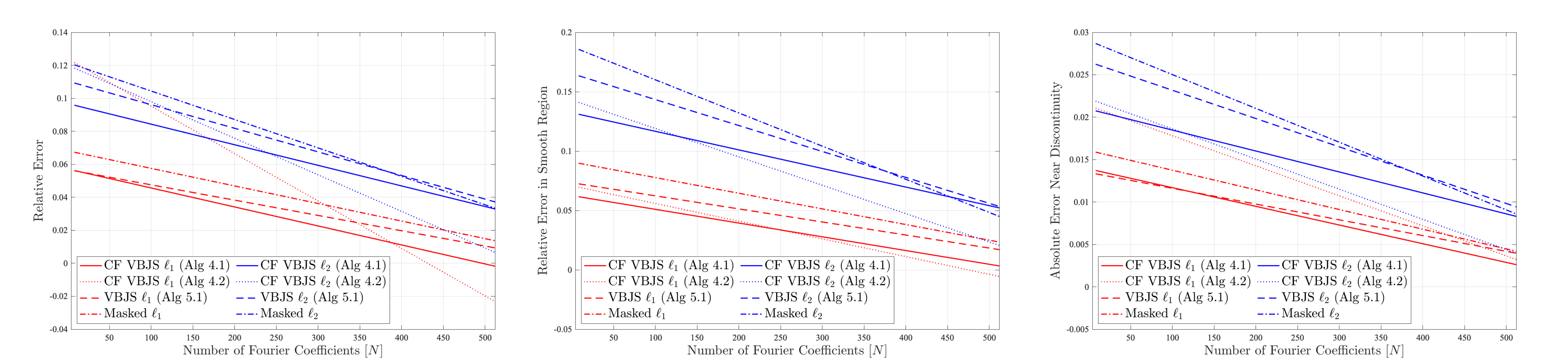}
\caption{Error comparison for reconstructing (\ref{eq:sawtooth}) using Algorithms \ref{alg:CF_VBJS_MMV}, \ref{alg:CF_VBJS_SMV}, and \ref{alg:VBJS_orig} along with the masking technique described in (\ref{eq:mask}), each with $p = 1$ and $2$, {and $\tau = 1$}.  The parameters used in (\ref{eq:four_coef}) and (\ref{eq:grid1D}) are $2N = N_x = 2^n$ with $n = 3,4,\dots,10$. The PA transform order in (\ref{eq:PA}) is $m = 2$, and noise level of $\text{SNR} = 5$ dB. (left) The relative error (\ref{eq:rel_error}) across the entire {domain}, (middle) the relative error (\ref{eq:rel_error}) in  a {smooth region, $x \in [-\pi,-1]$}, and (right) the absolute error (\ref{eq:abs_error}) at {$x_* \approx -0.1 $.} }
\label{fig:conv_analysis}
\end{figure}

For both the single and multi-measurement vector cases  we reconstruct the ramp function (\ref{eq:sawtooth}) from {$2N+1 = 2^n + 1$, $n = 3,4,\dots,10$  Fourier coefficients with a fixed noise level at $\text{SNR} = 5$ dB.  The spatial grid is chosen as  $N_x =  2N$.}  We consider solving the final reconstruction problem (\ref{eq:VBJS}) with $p=1$ and $p=2$ for each method. In every case, we choose $m =2$ for the order of the PA operator (\ref{eq:PA}), {as this is the ideal choice for piecewise linear functions}.  The accuracy of each reconstruction is calculated according to (\ref{eq:rel_error}) and (\ref{eq:abs_error}) over 50 independent trials, and the results are recorded and displayed in Figure \ref{fig:conv_analysis}.  Observe that for each specific algorithm, $p=1$ yields better results than  $p=2$, and using {\em any} algorithm with $p = 1$ is generally better than with $p = 2$ if  $N_x\geq 64$.  In each case, the masking technique (\ref{eq:mask}) yields the slowest convergence, emphasizing the necessity of the large separation of scales present in the proposed weighting scheme (\ref{eq:weights}). 

\subsubsection{{Robustness} experiment}

We now test the performance of the CF-VBJS technique as the SNR in (\ref{eq:SNR}) is decreased from $10$ dB to $-10$ dB, in $.1$ dB increments. Figure \ref{fig:SNR_analysis} displays the relative error (\ref{eq:rel_error}) in reconstructing (\ref{eq:sawtooth}) from a \textit{single} measurement of $2N + 1$ noisy Fourier coefficients (\ref{eq:four_coef}), with $2N = N_x = 128$ and compares our new approach to the original VBJS in Algorithm \ref{alg:VBJS_orig}, and the masked technique, that is using (\ref{eq:mask}) in (\ref{eq:VBJS}). The experiment was performed $50$ times for each of these techniques with $p = 1$ and $p = 2$.  Figure \ref{fig:SNR_analysis} displays the average relative relative error for each SNR level.  {Since we have a single measurement vector,  we use Algorithm \ref{alg:CF_VBJS_SMV} and again compute $J =10$ unique edge vectors in (\ref{eq:jump_approx_smv}) using exponential CFs defined by (\ref{eq:exp_cf}) with order $\alpha = 2j$ for $j =1,\dots,10$.  Observe that for both $p = 1$ and $2$, the CF VJBS method is more accurate than the other considered algorithms. It is also more robust for decreasing SNR.}

\begin{figure}
\centering
\includegraphics[width =\textwidth ]{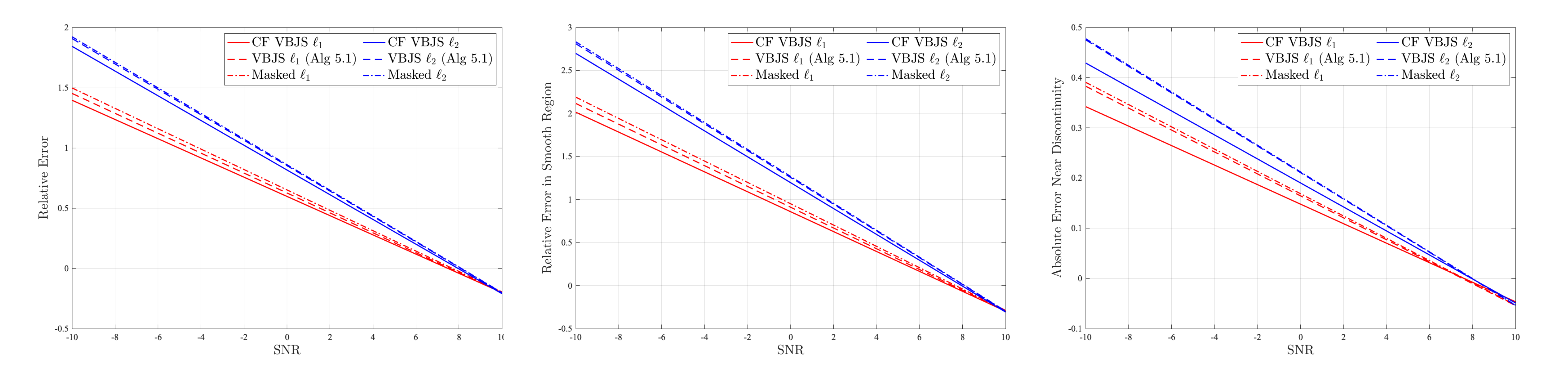}
\caption{Error in the reconstruction of (\ref{eq:sawtooth}) as a function of $\text{SNR}$ in (\ref{eq:SNR}). Methods compared include Algorithms \ref{alg:CF_VBJS_SMV} and \ref{alg:VBJS_orig} as well as the masking technique in (\ref{eq:mask}), each for $p = 1$ (red) and $2$ (blue).  We chose parameters $2N = N_x = 128$ corresponding to (\ref{eq:four_coef}) and (\ref{eq:grid1D}) respectively.  We also used PA transform order $m = 2$ in (\ref{eq:PA}). (left) The relative error (\ref{eq:rel_error}) across the entire reconstruction, (middle) the relative error (\ref{eq:rel_error}) for {$x \in [-\pi,-1]$}, and (right) the absolute error (\ref{eq:abs_error}) {at $x_* \approx -0.1$}.} 
\label{fig:SNR_analysis}
\end{figure}

\subsubsection{Image recovery in two dimensions}

The CF VBJS technique described in Algorithm \ref{alg:CF_VBJS_MMV} can be readily adapted for multi-dimensional problems by separately calculating edge maps in each direction.\footnote{We note that a two dimensional approach to edge detection via the concentration factor method was discussed in \cite{MGG,adcock2019joint}, but for uniform Fourier data it is more straightforward and efficient to use the dimension by dimension approach.}  Analogous to (\ref{eq:four_coef}), we now acquire $J$ measurements of Fourier data, $\hat{\bm  f} \in \mathbb{C}^{[2N+1]^2}$  given elementwise as
\begin{equation}
\label{eq:2D_FC}
\hat{f}_{k_x,k_y} = \frac{1}{4}\int_{-1}^{1}\int_{-1}^{1}f(x,y) e^{-i(k_x \pi x + k_y \pi y)} dx dy + \nu_{k_x,k_y}, \quad k_x,k_y = -N,...,N,
\end{equation}
where $\bm\nu \sim {\mathcal{CN}}[\mu,\sigma]$ is complex additive Gaussian noise with mean $\mu$ and standard deviation $\sigma$. 
We seek to recover  $f(x,y)$ on $[-1,1]^2$ on equally spaced grid points given by
\begin{eqnarray}
\label{eq:2Dgrid}
x_{j_x} &=& -1 + \frac{2(j_x-1)}{N_x}\nonumber\\
y_{j_y} &=& -1 + \frac{2(j_y-1)}{N_y},  
\end{eqnarray}
}
with $j_x = 1, \dots, N_x$, $j_y = 1,\dots, N_y$, and $N_x = N_y = 2N$.  Once again we apply the discrete Fourier transform (\ref{eq:fourmodel}) for the numerical model, where in this case we have the two-dimensional operator given by }
\begin{equation}
\label{eq:discrete_coef2D}
\mathcal{F}(k,j) =  \frac{e^{-i(k_x x_{j_x} + k_y y_{j_y})}}{N_xN_y},
\end{equation}
where $ k_x,k_y = -N,\dots,N$, $j_x = 1,\dots,N_x$, and $j_y = 1,\dots,N_y$. Algorithm \ref{alg:2D} describes the two-dimensional CF VBJS technique.
\begin{algorithm}[h!]
\caption{CF VBJS from Multiple Measurement Vectors in Two Dimensions.}
\label{alg:2D}
\begin{algorithmic}[1]
\STATE Acquire multiple measurement vectors $\hat{\bm f}^j$, $j = 1,\ldots,J,$ according to (\ref{eq:2D_FC}).
\STATE Using a single concentration factor $\sigma$, calculate $J$ jump function approximations in each direction as
\begin{align*}
\label{eq:jump_approx_2D}
\tilde{\bm g}^j_x &= i \sum_{|k_x|\leq N}\sum_{|k_y|\leq N} \hat{ f}^j_{k_x,k_y} {sgn}(k_x)\sigma\left(\frac{|k_x|}{N}\right)e^{i(k_x x + k_y y)}, \quad j = 1, \dots, J. \\
\tilde{\bm g}^j_y &= i \sum_{|k_x|\leq N}\sum_{|k_y|\leq N} \hat{ f}^j_{k_x,k_y} {sgn}(k_y)\sigma\left(\frac{|k_y|}{N}\right)e^{i(k_x x + k_y y)}, \quad j = 1, \dots, ,J.
\end{align*}
\STATE Use (\ref{eq:JS_matrix}) - (\ref{eq:weights}) to calculate the spatially adaptive weighting vector in each dimension as $w_x$ and $w_y$, respectively. The weighting vectors are reshaped into the appropriate $N_x \times N_y$ matrices, where $N_x$ and $N_y$ define the size of the spatial grid in (\ref{eq:2Dgrid}).
\STATE The final weighting matrix is then calculated as 
\begin{equation}
\label{eq:2Dweights}
w_{i,j} = \min_{i,j} \left[w_x(i,j),w_y(i,j)\right].
\end{equation}
\STATE Select the optimal measurement vector $\hat{\bm f}=\hat{\bm f}^{j^*}$ that solves (\ref{eq:jsdata_index}). 
\STATE With a choice of $p =1$ or $p = 2$ and PA transform order $m$, solve the CF VBJS reconstruction problem 
\begin{equation}
\label{eq:VBJS2D}
\bm{f}^* = \argmin{\bm q\in\mathbb{R}^{N_x\times N_y}} \left\{\frac{1}{p}||W\mathcal{L}^m \bm q ||_{p}^p + \frac{1}{2} ||\mathcal{F}\bm q - \hat{\bm{{f}}}||_2^2 \right\},
\end{equation}
where ${\mathcal F}$ is given in (\ref{eq:discrete_coef2D}).
\end{algorithmic}
\end{algorithm}

We use Algorithm \ref{alg:2D} to reconstruct the two-dimensional function $f(x,y)$ on $[-1,1]^2$ defined by 
\begin{equation}
\label{eq:2D_funct}
f(x,y) = \begin{cases} 
10\cos\left(\frac{3\pi}{2}\sqrt{x^2+y^2}\right), \quad \sqrt{x^2+y^2} \leq \frac{1}{2} \\
10\cos\left(\frac{\pi}{2}\sqrt{x^2+y^2}\right),  \quad \sqrt{x^2+y^2} > \frac{1}{2},
\end{cases}  
\end{equation}
with $2N = N_x = N_y = 128$.  The weighted Alternating Directions Method of Multipliers (ADMM) optimization algorithm described in \cite{gelb2018reducing} was used to solve (\ref{eq:VBJS2D}). 

\begin{figure}[htb]
\centering
\includegraphics[width =.8\textwidth ]{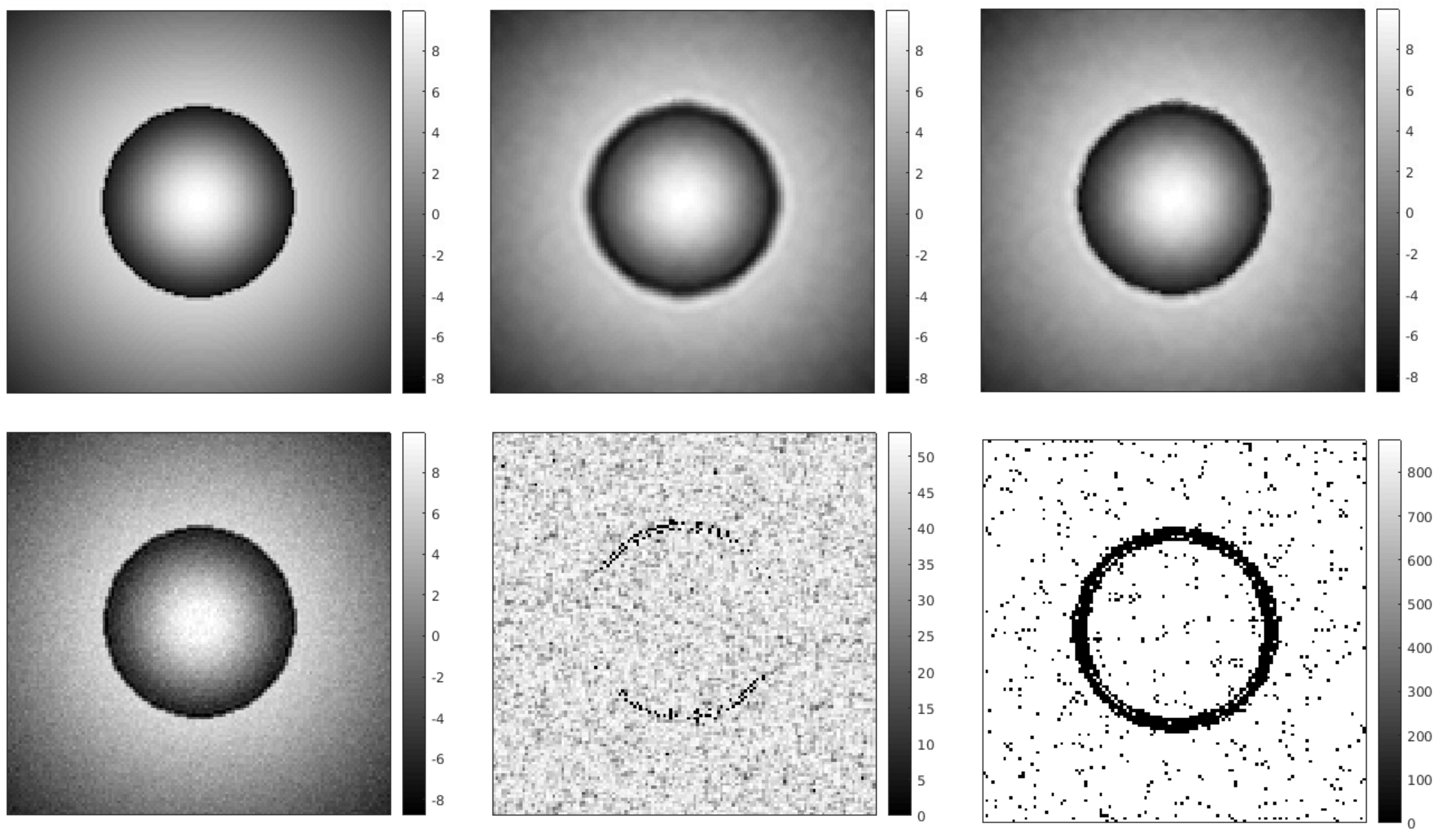}
\caption{{(top-left) The original, two-dimensional image. (bottom-left) The result of reconstructing from noisy Fourier data using the {filtered} inverse FFT. (top-middle) The result of reconstructing from noisy Fourier data using the VBJS technique described in \cite{gelb2018reducing}. (bottom-middle) The weighting matrix calculated using the technique described in \cite{gelb2018reducing}. (top-right) The result of reconstructing from noisy Fourier data using the CF VBJS technique with $p = 1$ and $m = 2$ in Algorithm \ref{alg:2D}. (bottom-right) The weighting matrix resulting from calculating (\ref{eq:2Dweights}).}}
\label{fig:2D}
\end{figure}

\begin{figure}[htb]
\centering
\includegraphics[width = .6\textwidth]{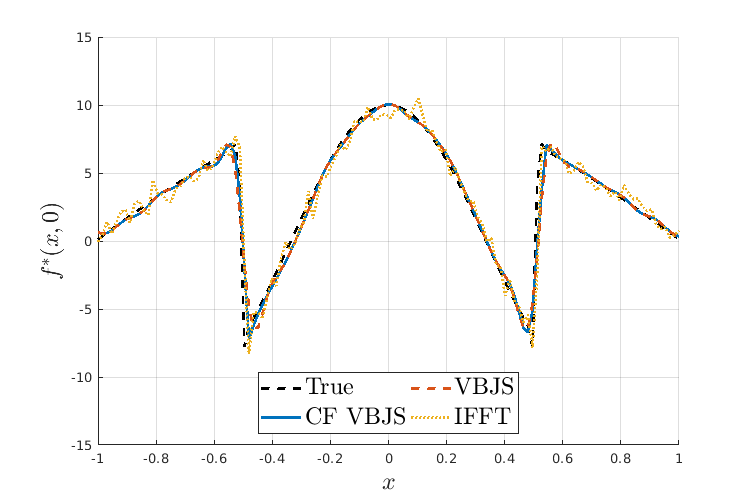}
\caption{{The one-dimensional cross-section of the results displayed in Figure \ref{fig:2D}. Cross sections are taken at $y = 0$.} }
\label{fig:2D_CS}
\end{figure}

{Figure \ref{fig:2D}(top-left) shows the underlying image of interest (\ref{eq:2D_funct}) from which we are given noisy Fourier data (\ref{eq:2D_FC}) with $\nu_{k_x,k_y}$ chosen so that $\text{SNR} = -10$ dB. The {filtered} IFFT is used to reconstruct the bottom-left image. Figure \ref{fig:2D} shows the reconstruction (top-middle) and weights (bottom-middle) calculated according to the VBJS technique described in \cite{gelb2018reducing}. Figure \ref{fig:2D}(bottom-right) shows the calculated weighting matrix (\ref{eq:2Dweights}), while Figure \ref{fig:2D}(top-right) displays the CF VBJS reconstruction with $p =1$ and $m = 2$ in (\ref{eq:VBJS}). As in the one-dimensional case, the CF VBJS technique quickly and accurately reconstructs the underlying image. This is further emphasized in Figure \ref{fig:2D_CS}, where we display the one-dimensional cross section of each reconstruction at $y=0$. The corresponding relative error (\ref{eq:rel_error}) and error near the a discontinuity (\ref{eq:abs_error}) are calculated in Table \ref{tab:2D_error}. In each case, the CF VBJS technique improves the accuracy of the reconstruction. 
}

\begin{table}[htb]
\begin{center}
\begin{tabular}{|c|c|c|c|}
\hline
 		 & Relative Error  & Relative Error & Absolute Error \\ 
 		 & (total) &  (in smooth regions)  &  ({at $x_*   \approx -0.57$})  \\ \hline
CF VBJS & 0.1534 & 0.0257  & 0.0167 \\ \hline
VBJS & 0.1769 & 0.0258  & 0.0366  \\ \hline 
IFFT & 0.1964 & 0.0725 & 0.1091  \\  \hline
\end{tabular}
\end{center}
\caption{{The error in each one-dimensional cross section (Figure \ref{fig:2D}) which results from estimating (\ref{eq:2D_funct}) using the the CF VBJS method described in Algorithm \ref{alg:2D}, the two-dimensional VBJS method from \cite{gelb2018reducing}, and the IFFT. We calculate both the relative error in the entire approximation and only in smooth regions, $|x|<0.2$. Also calculated is the absolute error (\ref{eq:abs_error}) at a point near an edge, $x_* \approx -0.57$}.}
\label{tab:2D_error}
\end{table}

\subsection{{Recovering signals from Fourier data with missing bands}}

We now seek to further analyze the iterative concentration factor (ICF) design developed in Section \ref{sec:ICF}. As reflected in Table \ref{tab:it_vs_exp_compare}, it is better to use the ICF approach to estimate the edge vectors (\ref{eq:approx_jump_funct}) for the VBJS approximation in (\ref{eq:VBJS}) when not all $2N+1$ Fourier coefficients (\ref{eq:four_coef}) are available; {For example, when certain frequency bands are inaccessible}. {To additionally test this result, we again fix $N = 64$ and gather $J = 4$ measurement vectors, where each measurement is missing Fourier coefficients within a band having bandwidth $b$. We repeat the experiments eight independent times,  increasing the bandwidth of the missing bands from $b = 2$ to $b = 16$ in two step increments. In each case we ensure that the missing bands are equally spaced throughout the $2N +1$ coefficients and that {$k = -N,0,N$ are always included in the set of known coefficients.} We compare edge vector approximations {given by (i) (\ref{eq:EdgeDetector}), which is used in the original  VBJS technique in {\cite{gelb2018reducing}}, (ii) the CF VBJS technique described in Algorithm \ref{alg:CF_VBJS_MMV} with exponential concentration factors (\ref{eq:exp_cf}), and (iii) the ICF VBJS technique outlined in Section \ref{sec:ICF}. We perform the reconstructions with $p = 1$ and $p=2$ in (\ref{eq:VBJS}). In each case, at every grid point, we calculate the log of the pointwise error as 
\begin{equation}
\label{eq:pw}
{\bm E}_{log} = \log_{10}|\bm{f}^* -\bm{f}|
\end{equation}
and display the results in Figure \ref{fig:pw_ICF}. It is evident that the ICF VBJS method produces more accurate results, especially in jump regions. Additionally, with the ICF VBJS method, we are able to exactly recover the unknown function in some cases.\footnote{{Again, this is due to the fact that the ICF method as described in (\ref{eq:ICF}) uses the ramp function as a template.}}   

\begin{figure}[htb]
\centering 
\includegraphics[width = \textwidth]{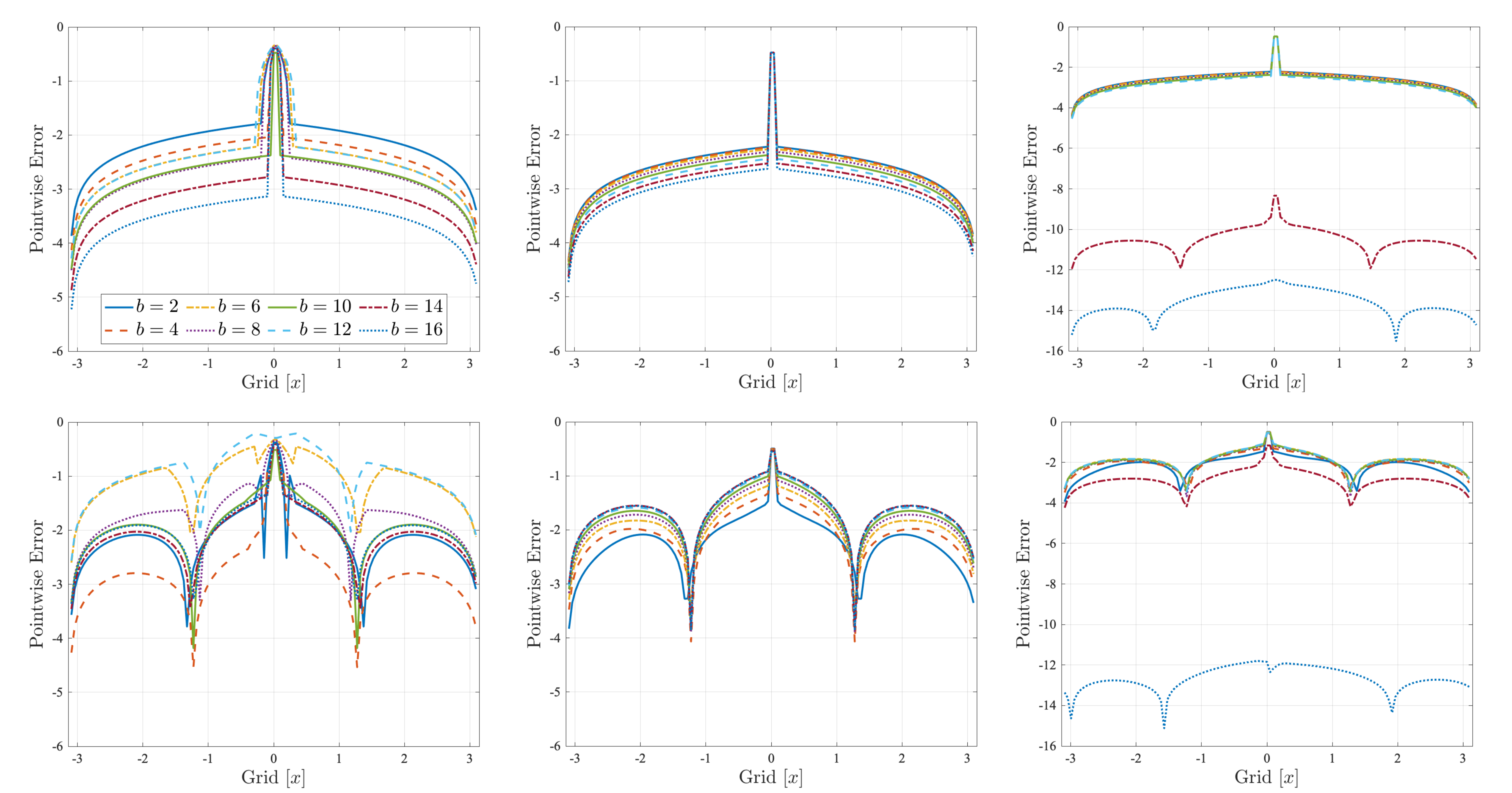}
\caption{{Pointwise error (\ref{eq:pw}) associated with reconstructing (\ref{eq:sawtooth}) when some bands of Fourier data are missing with bandwidth $b$.  Here we fix the number of Fourier coefficents to be $2N+1 = 129$. {The color legend is provided in the top-left figure}. (top) {$p = 1$} and (bottom) {$p=2$}.  The edge vectors were computed using (left) the orignal approach  proposed in \cite{gelb2018reducing, scarnati2018variance}, (middle) the CF VBJS method described in Algorithm \ref{alg:CF_VBJS_MMV}, and (right) the ICF VBJS technique explained in Section \ref{sec:ICF}. }}
\label{fig:pw_ICF}
\end{figure}

As a culminating experiment, we explore how the ICF VBJS method performs when a random percentage of bands are removed from each measurement. That is, we fix $N=64$ in (\ref{eq:four_coef}), and for 
$$\gamma = 0.05 +  l0.05, \quad l =0,...,18,$$
we randomly eliminate $\gamma (2N+1)$ bands of data. Note that $\gamma\in[0.05,0.95]$ represents the fraction of bands removed from the data. In each case, we reconstruct (\ref{eq:sawtooth}) with CF VBJS and ICF VBJS using $p=1$ and $p=2$ in (\ref{eq:VBJS}). {The relative reconstruction error (\ref{eq:rel_error}) and absolute error (\ref{eq:abs_error}) near the discontinuity in each case is reported in Figure \ref{fig:perc_bands}.} This experiment confirms that as expected, when bandwidths of data are removed, using $p=1$ provides greater accuracy than $p =2$ in (\ref{eq:VBJS}).  The ICF technique described in Section \ref{sec:ICF} furthermore improves the results.

\begin{figure}
\centering
\includegraphics[width=\textwidth]{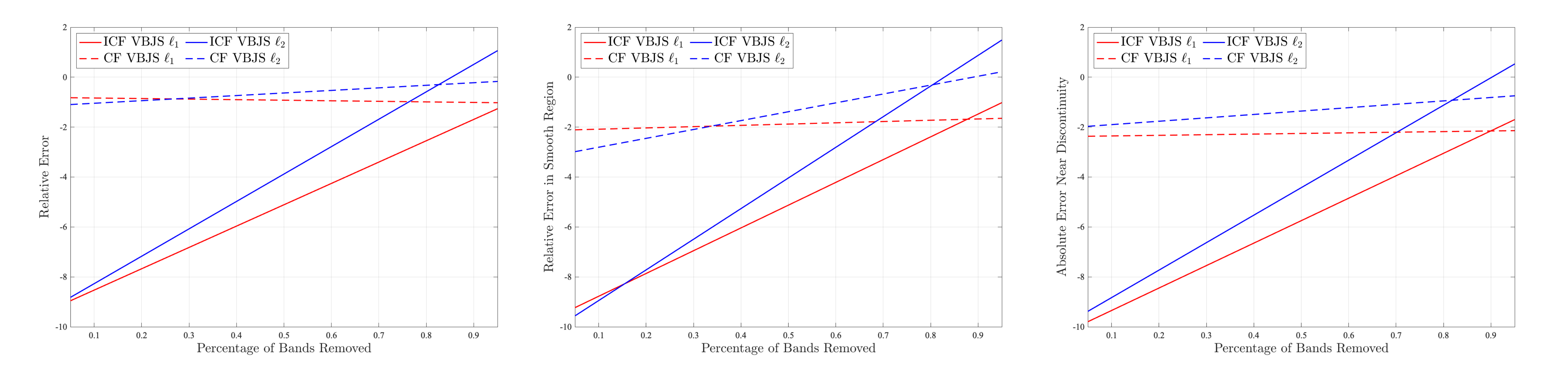}
\caption{{Relative and absolute error comparisons for the ICF and CF VBJS algorithms when a random fraction $\gamma$ of Fourier coefficients are removed from the data. (left) The relative error (\ref{eq:rel_error}) {over the total domain}, (middle) the relative error (\ref{eq:rel_error}) in a smooth region, {$|x| \leq 1$}, and (right) the absolute error (\ref{eq:abs_error}) {at $x_* \approx -0.1$.} Note that the y-axis of each plot has been log-scaled.}} 
\label{fig:perc_bands}
\end{figure}

%% file: conclusions.tex
This investigation introduced an accelerated version of the variance based joint sparsity (VBJS) recovery of images from Fourier data.  Given $J > 1$ data measurement vectors, the original version of VBJS provided in \cite{gelb2018reducing,adcock2019joint} reconstructed an underlying function $f$ by {\em first} producing $J$ individual reconstructions of $f$, typically using compressive sensing techniques.  This was followed by constructing a weighting vector for the regularization term based on the {\em joint sparsity} of the multi-measurements (in the sparsity domain).  Here, instead, we employed the concentration factor (CF) method on the given Fourier measurements to reconstruct $J$ edge vectors.  Already this yields two major advantages: (i) The CF edge detection method works directly on the Fourier data to produce edges and does not require solving an inverse problem, and (ii) no additional information is lost by having to first construct $J$ approximations of $f$.  We additionally enhanced our algorithm by employing the $minmod$ algorithm, which helps to reduce oscillatory repsonses in the CF approximations.  Consequently, we are able to produce a spatially variant weight vector that promotes the separation of scales apparent in the underlying function.  We note that the technique works when a single measurement vector (SMV) is acquired.  In this case $J$ approximations are obtained by using $J$ different concentration factors.  In addition to being more accurate, for both the SMV and MMV scenarios our new method is more efficient, since fewer inverse problems must be solved.

Our technique is also applicable when bands of Fourier samples may be missing from the acquired data.  In this case we use the iterative concentration factor (ICF) approach described in (\ref{eq:ICF}).  As the concentration factors can be seen as a weighting vector on the Fourier data, this process ensures that we are applying zero weights to the missing bands of Fourier data.  While dramatically improving the accuracy compared to the case when using the standard CFs, the cost is comparable to the original version of the VBJS in \cite{gelb2018reducing,adcock2019joint} since determining the CFs requires an iterative process.

There are several applications for which  accelerated VBJS technique will be useful.  For example, multiple measurement vectors of Fourier data are collected in synthetic aperture radar (SAR),  \cite{cetin2005sar} and in parallel magnetic resonance imaging (pMRI), \cite{ChunAdcockTalavagepMRI,UnserEtAlRealistic}.  Moreover, the accelerated VBJS can be helpful in employing change detection algorithms, since the difference in a reference and changed state can be more accurately measured.  Preliminary results demonstrate that this is indeed the case even in low SNR environments.  Future investigations will employ the accelerated VBJS technique in these applications.

%% file: appendixPA.tex
\label{appendix:PA}
While a variety of techniques can be used to approximate ${\bm g}$ in (\ref{eq:jump_funct_disc}) to be used in (\ref{eq:l1_reg}), for this investigation we employ the polynomial annihilation (PA) operator, ${\mathcal L}^m$, which is based on the polynomial annihilation edge detection method, \cite{archibald2005polynomial}.  The primary advantage in using the PA operator is that it is constructed to yield higher order approximations whenever $m > 1$. (When $m = 1$, it is equivalent to TV.)  This means that the approximation of ${\bm g}$ is more likely to be sparse, and specifically yields ${\bm g}_j \approx 0$ whenever the corresponding grid point $x_j$ falls in a smooth region of $f(x)$.  This is especially true when $f$ has more variation in smooth regions or when $\hat{\bm f}$ in (\ref{eq:four_coef}) is sparsely sampled.  Finally, as we will briefly describe in what follows, it is also easy to generate ${\mathcal L}^m$.  More information can be found in \cite{archibald2016image, stefan2010improved}.  

For ease of presentation we consider the one-dimensional case, where $f(x)$ is a piecewise smooth function on $[-\pi,\pi]$ (the finite domain is arbitrary) given on grid point values $x_j$, $j = 1,\cdots,N_x$.  They need not be uniform, although in our investigation we only consider uniform grid points given by (\ref{eq:grid1D}).  We also note that the the polynomial annihilation edge detection method can be described for any finite-dimensional function, \cite{archibald2005polynomial}. However, since images are usually defined on a Cartesian grid,  for  the purposes of defining a sparsifying transform operator to be used in (\ref{eq:l1_reg}), it was demonstrated in \cite{archibald2016image} that applying the PA transform operator dimension by dimension was accurate and efficient.  

We begin by defining the polynomial annihilation edge detection approximation as
\begin{equation}
\label{eq:EdgeDetector}
L^mf(x)=\frac{1}{q^{m}(x)}\sum_{x_j\in S_x}c_j(x)f(x_j),
\end{equation}
where $S_x$ is the local set of $m+1$ grid points from the set of given grid points about $x$, $c_j(x)$ are the polynomial annihilation edge detection coefficients, (\ref{eq:EdgeCoeffs}), and $q^m(x)$ is the normalization factor, (\ref{eq:normalization}).  Each parameter of the method can be further described as:
\begin{itemize}
\item $S_x$: For any particular cell $I_j=[x_j,x_{j+1})$, there are $m$ possible stencils, $S_x$ of size $m+1$, that contain the interval $I_j$.  For simplicity, we assume that the stencils are centered around the interval of interest, $I_j$, and are given by
$$S_{I_j} = \{x_{j - \frac{m}{2}},\cdots,x_{j+\frac{m}{2}}\},\hspace{.2in}
S_{I_j} = \{x_{j - \frac{m+1}{2}},\cdots,x_{j+\frac{m-1}{2}}\}$$
for $m$ even and odd respectively.  For non-periodic solutions the stencils are adapted to be more one sided as the boundaries of the interval are approached, \cite{archibald2005polynomial}.  To avoid cumbersome notation, we write $S_x$ as the generic stencil unless further clarification is needed.
\item $c_j(x)$:  The polynomial annihilation edge detection  coefficients, $c_j(x), \; j=1,\dots, m+1$, are constructed to annihilate polynomials up to degree $m$. They are obtained by solving the system
\begin{equation}
\label{eq:EdgeCoeffs}
\sum_{x_j\in S_x}c_j(x)p_\ell(x_j)= p_\ell^{(m)}(x), \; j=1,\dots, m+1,
\end{equation}
where $p_\ell$, $\ell=0,\dots,m,$ is a basis of for the space of polynomials of degree $\le m$.
\item $q^{m}(x)$:  The normalization factor, $q^m(x)$, normalizes the approximation to assure the proper convergence of $L^mf$ to the jump value at each discontinuity.  It is computed as
\begin{equation}
\label{eq:normalization}
q^m(x) = \sum_{x_j\in S_x^+} c_j(x),
\end{equation}
where $S_x^+$ is the set of points $x_j \in S_x$ such that $x_j \ge x$.
\end{itemize}

If the solution vector ${\bm f}$ is on uniform points, as it is in our case where $\{x_j\}_{j = 1}^{N_x}$ are defined in (\ref{eq:grid1D}), then there is an explicit formula for the polynomial annihilation edge detection coefficients, independent of location $x$, computed as (\cite{archibald2005polynomial})
\begin{equation}
\label{eq:UniformCoefficients}
c_j=\frac{m!}{\prod_{k=1, k \neq j}^{m+1}(j-k)\Delta x},        \;\; j=1,\dots ,m+1.
\end{equation}
Here $\Delta x = \frac{2\pi}{N_x}$.  We can now define the polynomial annihilation (PA) transform matrix, ${\mathcal L}^m$, as
\begin{equation}
\label{eq:PAtransform}
{\mathcal L}^m_{j,l}  = \frac{c(j,l)}{q^m(x_l)}, \hspace{.2in} 0 \le l \le N_x, \hspace{.2in} 0 \le j < N_x,
\end{equation}
where
$$c(j,l) = \left\{
     \begin{array}{lr}
       c_{j-l-\lfloor \frac{m}{2}\rfloor}, & 0<j-l-\lfloor \frac{m}{2}\rfloor+s(j,l)\leq m+1\\
       0 & \textrm{otherwise}
     \end{array}
   \right.
$$
and
$$s(j,l) = \left\{
     \begin{array}{lr}
       l-\lfloor \frac{m}{2}\rfloor, & \textrm{$l\leq\lfloor \frac{m}{2}\rfloor$}\\
        l+m-\lfloor \frac{m}{2}\rfloor-N_x, & \textrm{$l+m-\lfloor \frac{m}{2}\rfloor>N_x$}\\
       0 & \textrm{otherwise.}
     \end{array}
   \right.
$$
If the underlying image is known to be periodic, or is zero padded at the boundaries, a centered stencil can be used throughout the domain producing a circulant matrix ${\mathcal L}^m$.  A reduction of accuracy is expected near the boundaries for the non-periodic case due to the one-sided stencils.  

For example, assuming periodicity, the banded matrix ${\mathcal L}^1$  and ${\mathcal L}^3$ are given by
       \begin{equation}
        \label{eq:PAexamples}
                \mathcal{L}^1 = \begin{bmatrix}
                  1 & -1 &        &         &      \\
                    & 1  & -1     &         &      \\
                    &    & \ddots & \ddots  &       \\
                    &    &        &  1      & -1    \\
                 -1 &    &        &         & 1
                 \end{bmatrix}, \quad
                \mathcal{L}^3 = \frac{1}{2}\begin{bmatrix}
                  3  & -3     & 1      &        &                 & -1 \\
                  -1 &  3     & -3     & 1      &                 &     \\
                     & \ddots & \ddots & \ddots & \ddots          &     \\
                     &        &  1      &    3    &  -3           & 1   \\
                  1  &        &        & 1        &  3            & -3 \\
                 -3  & 1      &       &           &   -1          & 3
                 \end{bmatrix}.
        \end{equation}